# RINGS OF QUOTIENTS
# OF
# RINGS OF FUNCTIONS

N. J. FINE

L. GILLMAN

J. LAMBEK

$1.75

# Editors' preface

It has been forty years since these notes were printed (and forty-five since the results were announced). The results have been widely used, even though the original printing was not widely advertised and has long since been exhausted. Thus it seemed like a good idea to reprint it and make it freely available. Professors Gillman and Lambek were very supportive and readily assented to this. Unfortunately Professor Fine died in 1995 and, despite some effort, we have been unable to locate either his widow or his daughter. On the other hand, there is no copyright notice on the original edition and it would not seem that there is any possibility of commercial exploitation so that we have decided to go ahead nonetheless.

Since the original is a typescript, and not even high quality typescript, it was clear from the beginning that we would have it retyped by mathematicians who volunteered their services for the good of the community. We heartily thank these volunteers: Walter Burgess, Igor Khavkine, Chawne M. Kimber, Michelle L. Knox, Suzanne Larson, Ronald Levy, Warren Wm. McGovern, Jay Shapiro, and Eric R. Zenk. We also thank Gordon Mason and R. Grant Woods who have contributed to the project by proofreading the final result.

Michael Barr
Robert Raphael

Montreal, September, 2005

Added (2021-01-31) We would also like to thank Fred Dashiell for pointing out several typos and having provided the corrections. These are listed below.

- A new title page is a scan of the original.
- The fonts $\mathscr{Q}$ and Q are conformed to the original: $\mathscr{Q}$ for rings and Q (not $Q$) for spaces.
- The index more faithfully conforms to the original.
- A few typos in the 2005 edition and one typo in the original have been corrected. Any change to the original is indicated in a footnote.
- Science Citation Index now lists 91 articles which cite this document, of which 22 have been published in the last ten years.
- (Added 2023-03-01) The symbols $\cong$ (for isomorphic rings) and $\approx$ (for homeomorphic spaces) now conform to the original.

# Preface

These seminar notes contain the results announced in Abstracts 576-219, Notices Amer. Math. Soc. **7** (1960), 980; 61T-27, ibid. **8** (1961), 60; and 61T-28, ibid. **8** (1961), 61.

The rings of quotients recently introduced by Johnson and Utumi are applied to the ring $C(X)$ of all continuous real-valued functions on a completely regular space $X$. Let $Q(X)$ denote the maximal ring of quotients of $C(X)$; then $Q(X)$ may be realized as the ring of all continuous functions on the dense open sets in $X$ (modulo an obvious equivalence relation). In special cases (e.g., for metric $X$), $Q(X)$ reduces to the classical ring of quotients of $C(X)$ (formed with respect to the regular elements), but in general, the classical ring is only a proper subring of $Q(X)$.

A natural metric may be imposed on the ring $Q(X)$ (although $Q(X)$ does not, in general, become a topological ring). Let $\overline{\phantom{x}}$ denote the metric completion and let $*$ denote bounded functions. Then: $C(X)$ is a ring of quotients of $C^*(X)$, so that $Q(X) = Q(\beta X)$; $\overline{Q}^*(X)$ is a topological ring; $\overline{Q^*}(X) = \overline{Q}^*(X)$; and $\overline{Q}(X)$ is the maximal ring of quotients of $\overline{Q}^*(X)$. ($\beta X$ denotes the Stone-Čech compactification of $X$).

The ring $\overline{Q}(X)$ may be realized as the ring of all continuous functions on the dense $G_\delta$-sets in $\beta X$ (modulo an equivalence relation).

Let $K$ denote the maximal ideal space of $Q(X)$. The space $K$ is compact and extremally disconnected and is homeomorphic with the maximal ideal spaces of $Q^*(X)$, of $\overline{Q}(X)$, and of $\overline{Q}^*(X)$. In addition, $K$ is homeomorphic with the maximal ideal space of $E(Q(X))$, where $E(A)$ denotes the Boolean ring of all idempotents of a ring $A$. Incidentally, $E(\overline{Q}(X)) = E(Q(X))$; this ring is also isomorphic with the Boolean algebra $\mathfrak{B}(X)$ of all regular open subsets of $X$, so that $K$ is homeomorphic with the maximal ideal space of $\mathfrak{B}(X)$. (It follows that $K$ is the same for all separable metric spaces $X$ without isolated points.) Finally, $K$ may be realized as the inverse limit of the spaces $\beta U$, $U$ ranging over all dense opens sets in $X$; the spaces $\beta V$ with $V$ ranging over all dense open sets in $\beta X$; and the spaces $\beta S$ with $S$ ranging over all dense $G_\delta$'s in $\beta X$.

Results obtained in the course of the study include new proofs of known theorems, notably, the Stone-Nakano Theorem that $C(X)$ is conditionally complete as a lattice if and only if $X$ is extemally disconnected, and Artin's Theorem that in a formally real field, any element that is positive in every total order is a sum of squares.

# Acknowledgements

We began this work in 1960 while members of the Institute for Advanced Study. Fine was supported by Air Force Contract AF 18(603)-65 and, later, AF 49(638)-1007; Gillman by a National Science Foundation Senior Postdoctoral Fellowship and, later, by grants NSF G-14457 and NSF GP-3457; Lambek by a grant from the Canada Council.

We wish to thank several colleagues and students who read and criticized the manuscript, particularly Professor James Reid and Messrs. Kenneth Armstrong, Anthony Hager, and Donald Plank. We wish to thank Mrs. Ethel Brown, who typed the final version

The Pennsylvania State University
University of Rochester
McGill University

May 1965

# Contents

# 1. Rings of Quotients of Commutative Rings

**1.1.** MODULES. Let $A$ be a commutative ring with 1. Any ideal $I$ in $A$ may, of course, be regarded as an $A$-module. The set of all $A$-homomorphisms from $I$ into $A$ is denoted by $\operatorname{Hom}(I, A)$ or $\operatorname{Hom} I$.

The set $\operatorname{Hom} I$ is also an $A$-module. If $I'$ is an ideal, with $I' \subset I$, then the restriction map
$$\varphi \to \varphi | I' \qquad\qquad (\varphi \in \operatorname{Hom} I).$$
is a homomorphism of the module $\operatorname{Hom} I$ into the module $\operatorname{Hom} I'$.

Each $a \in A$ may be identified with the mapping $x \to ax$ $(x \in A)$; as such, it belongs to $\operatorname{Hom} A$.

**1.2.** THE CLASSICAL RING OF QUOTIENTS. The "classical" ring of quotients of $A$—to be denoted by $\mathfrak{Q}_{\mathrm{cl}}(A)$—is the ring of equivalence classes of formal quotients $c/d$ for $c \in A$ and $d$ a non-zero-divisor of $A$.

This construction has been generalized in the work of Johnson, Utumi, Findlay and Lambek [Johnson (1951), Johnson (1954), Johnson (1957), Utumi (1956), Findlay & Lambek (1958)]. Let $d$ be any non-zero-divisor of $A$. For each $c \in A$, a mapping $c/d$ may be defined on the principal ideal $(d)$ as follows:
$$\frac{c}{d}(dx) = cx;$$
then $c/d \in \operatorname{Hom}(d)$. Conversely, if $\varphi \in \operatorname{Hom}(d)$, then $\varphi(dx) = \varphi(d) \cdot x$, so that $\varphi$ acts like $\varphi(d)/d$. Thus one may construct $\mathfrak{Q}_{\mathrm{cl}}(A)$ by taking the direct limit of the modules $\operatorname{Hom}(d)$ and then defining a suitable ring structure on this direct limit. The generalization consists in considering modules $\operatorname{Hom} D$ for arbitrary "dense" ideals $D$, not just for the principal ideals generated by non-zero-divisors. There are further generalizations - e.g., for noncommutative rings. However, our discussion will be limited to commutative rings with unity.

**1.3.** DENSE IDEALS. An ideal in $A$ is said to be (*rationally*) *dense* if its only annihilator in $A$ is 0. Note that a principal ideal $(d)$ is dense precisely when $d$ is a non-zero-divisor.

More generally, we can speak of denseness of any subring. If $A$ is a subring of $B$, then we shall say that $A$ is (*rationally*) *dense in* $B$ provided that $A$ has no nonzero annihilator in $B$.

Any subring containing a dense subring is dense. The product of a dense ideal with a dense subring is dense. Hence the intersection of two dense ideals is dense.

The following lemma will be useful. An ideal is said to be *large* if it has nonzero intersection with every nonzero ideal—equivalently, with every nonzero principal ideal. Obviously, every dense ideal is large. The ring $A$ is *semi-prime* if it has no nilpotent ideal except (0)—equivalently (for commutative $A$), if it has no nilpotent element except 0.

LEMMA.    *A is semi-prime if and only if every large ideal is dense.*

PROOF. Assume $A$ semi-prime, and let $L$ be large and $a \neq 0$. Since $L \cap (a) \neq 0$, we have $0 \neq (L \cap (a))^2 \subset La$, which shows that $L$ is dense. Conversely, assume every large ideal dense, and let $I$ be an ideal such that $I^2 = 0$. Then for any $a$, we have $AI \subset I\hat{\ }$, where $I\hat{\ } = \{x \in A : xI = 0\}$, the annihilator of $I$. Hence if $(a) \cap I\hat{\ } = 0$, then $aI \subset (a) \cap I\hat{\ } = 0$—that is, $a \in I\hat{\ }$; consequently, $a \in (a) \cap I\hat{\ } = 0$. This shows that $I\hat{\ }$ is large. Therefore $I\hat{\ }$ is dense. Since $II\hat{\ } = 0$, this yields $I = 0$. ∎

**1.4.  RINGS OF QUOTIENTS; RATIONAL EXTENSIONS.** Let $B$ be a commutative ring containing $A$ and having the same unity element. For $b \in B$, we write

$$b^{-1}A = \{a \in A : ba \in A\}.$$

Obviously, $b^{-1}A$ is an ideal in $A$. For $b = 0$ or $1$, $b^{-1}A$ is dense: it is $A$.

DEFINITION.    $B \ (\supset A)$ *is a* ring of quotients *or* rational extension *of* $A$ *provided that for every* $b \in B$, $b^{-1}A$ *is dense in* $B$—*that is,*

(1)         *for* $0 \neq b' \in B$, *there exists* $a \in A$ *such that* $ba \in A$ *and* $b'a \neq 0$.

A ring without a proper rational extension is said to be *rationally complete*. (For a (commutative) semi-prime ring $A$, this coincides with *self-injective*, that is to say, injective as an $A$-module.)

Let $A \subset B \subset C$; if $C$ is a ring of quotients of $A$, then obviously $C$ is a ring of quotients of $B$ and $B$ is a ring of quotients of $A$. The converse is also true. To see this, we note first that if $B$ is a ring of quotients of $A$, then the element $a$ in (1) can always be chosen so that $b'a \in A$: given $a_1$ as in (1), choose $a_2 \in b'^{-1}A$ such that $b'a_1a_2 \neq 0$; then take $a = a_1a_2$. Suppose now that $C$ is a ring of quotients of $B$ and $B$ is a ring of quotients of $A$. Given $c \in C$ and $0 \neq c' \in C$, choose $b \in B$ such that $cb \in B$ and $0 \neq c'b \in B$; now choose $a_1 \in A$ such that $cba_1 \in A$ and $c'ba_1 \neq 0$, and then pick $a_2 \in A$ such that $ba_2 \in A$ and $c'ba_1a_2 \neq 0$; finally, take $a = ba_1a_2$. Then $ca \in A$ and $c'a \neq 0$, q.e.d.

Next, we observe: *if $B$ is a ring of quotients of $A$ and $D$ is a dense ideal in $A$, then $D$ is dense in $B$.* For consider any $b \neq 0$ in $B$. Choose $a \in A$ for which $0 \neq ba \in A$ and then $d \in D$ for which $(ba)d \neq 0$. Then $ad \in D$ and $b(ad) \neq 0$.

Finally, if $B$ is a ring of quotients of $A$ and $E$ is a dense ideal in $B$, then $E \cap A$ is dense in $B$. For let $b(E \cap A) = 0$. Then for each $e \in E$, we have $be(e^{-1}A) = 0$; since $e^{-1}A$ is dense, this yields $be = 0$. Thus, $bE = 0$ and therefore $b = 0$.

1.5  LEMMA    Let $B \supset A$.

(1) $B$ is a ring of quotients of $A$ if and only if for each nonzero $b \in B$, $b^{-1}A$ is dense in $A$ and $b \cdot (b^{-1}A) \neq 0$.

(2) If $b \cdot (b^{-1}A) \neq 0$ for all nonzero $b \in B$, then each ideal $b^{-1}A$ is large.

PROOF. (1). The necessity is trivial. Conversely, to show that $b^{-1}A$ is dense in $B$, consider any nonzero $b' \in B$. Pick $x \in b'^{-1}A$ for which $b'x \neq 0$; then pick $a \in b^{-1}A$ for which $b'xa \neq 0$.

(2). Let $a \in A$, with $a \neq 0$, be given. If $a \in b^{-1}A$, we are done. If $a \notin b^{-1}A$, then $ba \neq 0$, and so there exists $a' \in (ba)^{-1}A$ for which $baa' \neq 0$; then $0 \neq aa' \in (a) \cap b^{-1}A$. ∎

THEOREM.    Let $B \supset A$. If $A$ is semi-prime, then $B$ is a ring of quotients of $A$ if and only if $b \cdot (b^{-1}A) \neq 0$ for all nonzero $b \in B$—that is, for $0 \neq b \in B$, there exists $a \in A$ such that $0 \neq ba \in A$.

PROOF. This is an immediate consequence of the lemmas. ∎

**1.6.  CLOSED FAMILIES OF $A$-MODULES.**  We shall say that a family $\mathcal{D}$ of dense ideals in $A$ is *closed* provided that $A \in \mathcal{D}$ and the product of any two members of $\mathcal{D}$ is a member of $\mathcal{D}$. Thus, the smallest closed family is $\{A\}$. (Since $1 \in A$, $A \cdot A = A$.) The largest closed family is

$\mathcal{D}_0(A)$, the family of *all* dense ideals in $A$.

If $D$ and $D'$ are dense ideals, with $D \supset D'$, then the restriction homomorphism $\varphi \to \varphi|D'$ from $\operatorname{Hom} D$ into $\operatorname{Hom} D'$ is a *monomorphism*. For if $0 \neq \varphi \in \operatorname{Hom} D$, then $\varphi(d) \neq 0$ for some $d \in D$; since $D'$ is dense, there exists $d' \in D'$ such that $0 \neq \varphi(d) \cdot d' = \varphi(dd')$; therefore $\varphi|D' \neq 0$. We abuse notation to write $\operatorname{Hom} D \subset \operatorname{Hom} D'$.

Let $\mathcal{D}$ be a closed family of dense ideals. By a *standard* family of $A$-modules, we shall mean a family $(D^h)_{D \in \mathcal{D}}$ satisfying:

($\alpha$) $D^h$ is a submodule of $\operatorname{Hom} D$,

($\beta$) $1 \in A^h$

($\gamma$) if $D \supset D'$, then $D^h \subset D'^h$,

($\delta$) if $\varphi_1, \varphi_2 \in D^h$, then $\varphi_1 \circ \varphi_2 \in (DD)^h$.

Of course, $1 \in \operatorname{Hom} A$. If $D \supset D'$, then trivially $D^h \subset \operatorname{Hom} D'$. In ($\delta$), $\varphi_1 \circ \varphi_2$ is indeed defined on $DD$ (since $\varphi_2(dd') = d \cdot \varphi_2(d')$), and obviously it belongs to $\operatorname{Hom}(DD)$. In particular, then, $(\operatorname{Hom} D)_{D \in \mathcal{D}}$ is a standard family of $A$-modules.

**1.7. Direct limits.** We now consider the direct limit module

$$\mathcal{Q}_h(A) = \varinjlim_{D \in \mathcal{D}} D^h.$$

Because of ($\gamma$), $\mathcal{Q}_h(A)$ may be thought of as $\cup_{D \in \mathcal{D}} D^h$, where we identify $\varphi_1 \in D_1^h$ with $\varphi_2 \in D_2^h$ whenever $\varphi_1$ and $\varphi_2$ agree on $D_1 D_2$. The module operations in $\mathcal{Q}_h(A)$ then reduce to the operations within each $D^h$.

To obtain a ring structure on $\mathcal{Q}_h(A)$, we define

$$\varphi_1 \cdot \varphi_2 = \varphi_1 \circ \varphi_2,$$

which makes sense because of ($\delta$). Then $\mathcal{Q}_h(A)$ becomes a commutative ring with unity element $1 \in A$. By ($\beta$), $a = a \cdot 1 \in A^h$ for each $a \in A$; hence $A \subset \mathcal{Q}_h(A)$.

Let $(D^h)_{D \in \mathcal{D}}$ and $(E^k)_{E \in \mathcal{E}}$ be standard families such that $\mathcal{D} \subset \mathcal{E}$, and $D^h \subset D^k$ for $D \in \mathcal{D}$; then $\mathcal{Q}_h(A) \subset \mathcal{Q}_k(A)$. Evidently,

$$\mathcal{Q}(A) = \varinjlim_{D \in \mathcal{D}_0(A)} \operatorname{Hom} D.$$

is the largest of all the rings $\mathcal{Q}_h(A)$, and $A$ itself is the smallest.

**1.8. Representation theorem.** The rings of quotients of $A$ are precisely the direct limit rings $\mathcal{Q}_h(A)$. Specifically:

(1) If $(D^h)_{D \in \mathcal{D}}$ is any standard family of $A$-modules, then

$$\mathcal{Q}_h(A) = \varinjlim_{D \in \mathcal{D}} D^h$$

is a ring of quotients of $A$.

(2) If $B$ is any ring of quotients of $A$, then $(B \cap \operatorname{Hom} D)_{D \in \mathcal{D}_0(A)}$ is a standard family and

$$B = \varinjlim_{D \in \mathcal{D}_0(A)} (B \cap \operatorname{Hom} D).$$

PROOF. (1). Let $\varphi \in \mathcal{Q}_h(A)$; we are to prove that $\varphi^{-1}A$ is dense in $\mathcal{Q}_h(A)$. Given $0 \neq \varphi' \in \mathcal{Q}_h(A)$, pick $D$ such that $\varphi, \varphi' \in D^h$ and then choose $d \in D^h$ for which $\varphi'(d) \neq 0$. Since $(\varphi \circ d)(a) = \varphi(da) = \varphi(d) \cdot a$ for $a \in A$, we have $\varphi \cdot d = \varphi \circ d = \varphi(d) \in A$; hence $d \in \varphi^{-1}A$ and $\varphi' \cdot d = \varphi'(d) \neq 0$.

(2). This is a consequence of (1) (see next section), but the following proof also seems of interest. If $b \in B$, $D \in \mathcal{D}_0(A)$, and $b|D \in \operatorname{Hom} D$, then $b$ is uniquely determined by $b|D : b \neq 0$ implies $bD \neq 0$ ($D$ being dense in $B$), i.e., $b|D \neq 0$, so that the mapping $b \to b|D$ is one-one; hence, abusing notation to write $b \in \operatorname{Hom} D$, we see that $B \cap \operatorname{Hom} D$ has meaning. Next, it is obvious that $(B \cap \operatorname{Hom} D)_{D \in \mathcal{D}_0}$ is a standard family (1.6). Finally, each $b^{-1}A$ is dense and $b \in \operatorname{Hom} b^{-1}A$, which shows that $B \subset \cup_{D \in \mathcal{D}_0}(B \cap \operatorname{Hom} D)$ [$= \varinjlim_D (B \cap \operatorname{Hom} D)$]; and the reverse inclusion is trivial. ∎

**1.9 COROLLARY** *$A$ has a largest ring of quotients—namely,*

$$\mathcal{Q}(A) = \varinjlim_{D \in \mathcal{D}_0(A)} \operatorname{Hom} D.$$

This follows from (1) and the remark at the end of 1.7. The result implies (2): $(B \cap \operatorname{Hom} D)_{D \in \mathcal{D}_0}$ is a standard family, and $B \subset \mathcal{Q}(A)$.

The ring $\mathcal{Q}(A)$ is called the *maximal ring of quotients* or *rational completion* of $A$; evidently, $\mathcal{Q}(A)$ is rationally complete. There is also a smallest ring of quotients: $A$ itself.

**1.10 COROLLARY** *The classical ring of quotients—$\mathcal{Q}_{\mathrm{cl}}(A)$—is given by*

$$\mathcal{Q}_{\mathrm{cl}}(A) = \varinjlim_{(d) \in \mathcal{D}_0(A)} \operatorname{Hom}(d).$$

This follows from the discussion in 1.2.

**1.11.** REGULAR RINGS. The ring $A$ (commutative, with 1) is *regular* if for each element $a$, there exists an element $b$ (in general, not unique) such that $a^2 b = a$. Since $a(1 - ab) = 0$, every element is either a zero-divisor or a unit. If $a \neq 0$, then $1 - ab$ is not a unit and so belongs to some maximal ideal $M$, whence $a \notin M$; this shows that $A$ is *semi-simple* (i.e., the intersection of all maximal ideals is 0). Hence $A$ is semi-prime.

**THEOREM.** *If $A$ is semi-prime, then $\mathcal{Q}(A)$ is regular [Johnson (1951)].*

PROOF. Consider any $b \in \mathcal{Q}(A)$. Then $b \in \operatorname{Hom}(D, A)$ for some dense ideal $D$ in $A$. Let $K$ denote the kernel of $b$ and $K\hat{\ }$ its annihilator:

$$K = \{d \in D : bd = 0\}, \quad K\hat{\ } = \{a \in A : aK = 0\}.$$

Since $A$ is semi-prime, $K \cap K\hat{\ } = 0$; it follows that $b$ is one-one on $D \cap K\hat{\ }$. Put $E = b(D \cap K\hat{\ })$. Since $(E + E\hat{\ })\hat{\ } = E\hat{\ } \cap E\hat{\ }\hat{\ } = 0$, $E + E\hat{\ }$ is dense. (Similarly, $K + K\hat{\ }$ is dense.) Define $c \in \operatorname{Hom}(E + E\hat{\ }, A)$ by means of:

$$c(bd) = d \text{ for } d \in D \cap K\hat{\ }, \quad c(E\hat{\ }) = 0.$$

Then $b^2cd = bd$ for $d \in D \cap K\hat{}$; and the equation holds trivially for $d \in D \cap K$. Consequently, $b^2c - b$ annihilates the dense ideal $D \cap (K + K\hat{})$ and hence is 0. ∎

The proof of this special case of [Johnson (1951), Theorem 2] is of some interest in that it does not require Zorn's lemma.

# 2. Rings of Quotients of $C(X)$ and $C^*(X)$

**2.1.** RINGS OF FUNCTIONS. Let $X$ be a topological space. Let $C(X)$ denote the ring of all continuous functions from $X$ into the real line **R**, under the pointwise operations; the subring of bounded functions in $C(X)$ is denoted by $C^*(X)$. It is obvious that $C$ and $C^*$ are commutative rings with zero and unity elements the constant functions 0 and 1. They are also lattices under the pointwise definition of order.

It is obvious that $C(X)$ is semi-prime. (In fact, it is easy to see that $C(X)$ is semi-simple.) Hence Theorem 1.5 yields:

> A ring $B \supset C(X)$ is a ring of quotients of $C(X)$
> if and only if $b(b^{-1}C(X)) \neq 0$ for all nonzero $b \in B$.

In studying $C(X)$, we assume without loss of generality that the space $X$ is *completely regular*, i.e., it is a Hausdorff space such that for any neighborhood $U$ of a point $x$, some function in $C^*(X)$ vanishes outside $U$ but not at $x$; using lattice properties, it is then easy to construct a continuous function that vanishes on some *neighborhood* of $X - U$ and has a constant nonzero value on some neighborhood of $x$. It is also easy to prove that all subspaces of completely regular spaces are completely regular. (See [Gillman & Jerison (1960)].)

If $S$ is dense in $X$, then the homomorphism $f \mapsto f|S$ from $C(X)$ into $C(S)$ is a monomorphism; we abuse notation to write $C(X) \subset C(S)$.

The family of all *dense open* sets in $X$ is denoted by $\mathcal{V}_0(X)$.

**2.2.** ZERO-SETS AND COZERO-SETS. The zero-set of a continuous function $f$—i.e., its set of zeros—is denoted by $Z_X(f)$ or $Z(f)$. Obviously, every zero-set is closed: $Z(f) = f^{\leftarrow}(0)$, the preimage of a point in **R**. (The arrow indicates inverse of mapping.) In addition, every zero-set is a $G_\delta$ (i.e., a countable intersection of open sets) in $X$, since $\{0\}$ is a $G_\delta$ in **R**. The cozero set of $f$ is denoted by $\text{coz}_X f$ or $\text{coz } f$:

$$\text{coz } f = X - Z(f) = \{x \in X : f(x) \neq 0\}.$$

For any ideal $I$ in $C$ (or $C^*$), we write

$$Z(I) = \bigcap_{f \in I} Z(f), \qquad \text{coz } I = \bigcup_{f \in I} \text{coz } f.$$

Obviously, $\operatorname{coz} I$ is open (and, dually, $Z(I)$ is closed). Moreover, every open set $U$ is of the form $\operatorname{coz} I$ for some ideal $I$: define $I = \{f \in C(X) : \operatorname{coz} f \subset U\}$; then $I$ is an ideal and $\operatorname{coz} I \subset U$; by complete regularity, $\operatorname{coz} I = U$.

Evidently, $\operatorname{coz}(f)$ is the cozero-set $\operatorname{coz} f$. However, an arbitrary open set is not necessarily a cozero-set; so when an ideal $I$ is not principal, $\operatorname{coz} I$ need not be a cozero-set (in spite of the notation). Cf. 3.2.

THEOREM.    An ideal $D$ in $C$ (or $C^*$) is (rationally) dense if and only if $\operatorname{coz} D$ is (topologically) dense in $X$.

PROOF. The following are successively equivalent: $D$ is dense; for all $g \in C$, $gD = 0$ implies $g = 0$; for all $g \in C$, $Z(g) \supset \operatorname{coz} D$ implies $g = 0$; $\operatorname{coz} D$ is dense (by complete regularity). ∎

2.3    THEOREM

(1) $C(X)$ is a ring of quotients of $C^*(X)$.

(2) $C(X)$ and $C^*(X)$ have the same maximal ring of quotients and the same classical ring of quotients.

(3) If $V \in \mathcal{V}_0(X)$, then $C(V)$ is a ring of quotients of $C(X)$.

PROOF. If $0 \neq g \in C(X)$, then $f = 1/(1+g^2) \in C^*(X)$ and $0 \neq gf \in C^*(X)$. This yields (1), and the first half of (2) follows. For the second half of (2), we note that any quotient $g/h$ of two functions in $C(X)$ is equivalent (see 1.2) to the quotient $gf/hf$, where $f = 1/(1+g^2+h^2)$, of two functions in $C^*(X)$. (Cf. 3.8.)

For (3), we note first that $C(V) \supset C(X)$, because $V$ is dense. Now consider any $h \in C(V)$. For each $v$ in the open set $V$, there exists $f \in C^*(X)$ that vanishes on a neighborhood of $X - V$ but not at $v$; then $v \in \operatorname{coz} f$ and $hf \in C(X)$. In case $v \in \operatorname{coz} h$, we have $hf \neq 0$, which yields (3). In general, if we put $f' = f/(1+h^2f^2)$, then $f' \in C^*(X)$ and $hf' \in C^*(X)$, i.e., $f' \in h^{-1}C^*(X)$; and $v \in \operatorname{coz} f'$. So we have (for later use):

(4) $\qquad\qquad$ If $V \in \mathcal{V}_0(X)$ and $h \in C(V)$, then $V \subset \operatorname{coz} h^{-1}C^*(X)$. ∎

It follows from (3) that $C(V)$ and $C(X)$ have the same maximal ring of quotients. However, they do not in general have the same classical ring of quotients; see 3.7.

**2.4.    DIRECT LIMITS.** Let $\mathcal{S}$ be a family of nonvoid subsets of $X$. When $\mathcal{S}$ is a filter base (i.e., when $\mathcal{S}$ is closed under finite intersection), we are invited to consider the direct limit ring $\varinjlim_{S \in \mathcal{S}} C(S)$, with respect to the restriction homomorphisms $f \mapsto f|S'$, when $f \in C(S)$ and $S \supset S'$. When $\mathcal{S}$ is a family of *dense* sets, all these homomorphisms are one-one, and $\varinjlim_{S \in \mathcal{S}} C(S)$ may be thought of as $\bigcup_{S \in \mathcal{S}} C(S)$, where we identify $f_1 \in C(S_1)$ with $f_2 \in C(S_2)$ whenever $f_1$ and $f_2$ agree on $S_1 \cap S_2$.

For the time being, we shall confine our attention to the filter base $\mathcal{V}_0(X)$ of all dense open sets. If $h \in \varinjlim_{V \in \mathcal{V}_0} C(V)$, $h \neq 0$, then $h \in C(V)$ for some $V \in \mathcal{V}_0$; since $C(V)$ is a ring of quotients of $C(X)$, we have $h \cdot (h^{-1}C(X)) \neq 0$. This shows that $\varinjlim_{V \in \mathcal{V}_0} C(V)$, is also a ring of quotients of $C(X)$. Consequently, it is a subring of the maximal ring of quotients $\mathcal{Q}(C(X))$. Explicitly: if $h \in C(V)$, then $D = h^{-1}C(X)$ is a dense ideal in $C(X)$ and $h \in \operatorname{Hom} D$. Conversely, if $D$ is dense and $\varphi \in \operatorname{Hom} D$, then $\varphi \in C(V)$ for some $V$, as we now show.

**2.5 Lemma** *Let $A$ be a subring of $C(X)$; for any ideal $D$ in $A$, we have $\operatorname{Hom} D \subset C(\operatorname{coz} D)$.*

PROOF. Let $\varphi \in \operatorname{Hom} D$ be given. For $x \in \operatorname{coz} D$, choose $d \in D$ for which $d(x) \neq 0$, and define
$$g(x) = \frac{\varphi(d)(x)}{d(x)}.$$

Since $\varphi(d') \cdot d = \varphi(d) \cdot d'$, this definition is independent of $d$. For each $x \in \operatorname{coz} D$, $g$ agrees with a continuous function on a neighborhood (namely, $\operatorname{coz} d$) of $x$; therefore $g$ is continuous on its domain $\operatorname{coz} D$. Finally, consider any $d \in D$. For each $x \in \operatorname{coz} D$, we have, via suitable $d'$,
$$\varphi(d)(x) = \frac{\varphi(d')(x)}{d'(x)} d(x) = g(x) \cdot d(x);$$

hence $\varphi(d) = g \cdot d$ and therefore $\varphi = g$. ∎

REMARK. It is not true in general that $\operatorname{Hom} D = C(\operatorname{coz} D)$, even when $A = C(X)$ and $D$ is dense (although it is easily seen that $\operatorname{Hom} D \supset C^*(\operatorname{coz} D)$). Take $X = \mathbf{R}, A = C(\mathbf{R}), D = \{f \in C(\mathbf{R}) : f(0) = 0\}, d(r) = r, g(x) = 1/x^2$; then $gd \notin C(\mathbf{R})$, i.e., $g \notin \operatorname{Hom} D$.

**2.6. Q(X) AND Q_{cl}(X).** We write
$$\operatorname{Q}(X) = \mathcal{Q}(C(X)), \qquad \operatorname{Q}_{cl}(X) = \mathcal{Q}_{cl}(C(X)).$$

THEOREM. [REPRESENTATION THEOREM]

(1) $\operatorname{Q}(X) = \varinjlim_{V \in \mathcal{V}_0(X)} C(V)$.

(2) $\operatorname{Q}_{cl}(X) = \varinjlim_V C(V)$, $V$ ranging over all dense cozero-sets in $X$.

PROOF. Because $\operatorname{Hom} D \subset C(\operatorname{coz} D)$ for each $D$, we have

$$\mathrm{Q}(X) = \varinjlim_{D \in \mathcal{D}_0} \operatorname{Hom} D \subset \varinjlim_{D \in \mathcal{D}_0} C(\operatorname{coz} D) = \varinjlim_{V \in \mathcal{V}_0} C(V) \subset \mathrm{Q}(X).$$

This yields (1). To establish (2), we begin as above; we have

$$\mathrm{Q}_{\mathrm{cl}}(X) = \varinjlim_{(d) \in \mathcal{D}_0} \operatorname{Hom}(d) \subset \varinjlim_{(d) \in \mathcal{D}_0} C(\operatorname{coz} d) = \varinjlim C(V),$$

where $V$ ranges over the dense cozero-sets in $X$. We are to establish the reverse inclusion; thus, given $(d) \in \mathcal{D}_0$ and $f \in C(\operatorname{coz} d)$, we are to find $(d') \in \mathcal{D}_0$ such that $f \in \operatorname{Hom}(d')$. Define $d'$ to be equal to $d/(1 + f^2)$ on $\operatorname{coz} d$ and 0 otherwise; then $d' \in C(X)$. Since $\operatorname{coz} d' = \operatorname{coz} d$, a dense set, $(d') \in \mathcal{D}_0$. Finally, consider any member of $(d')$; it has the form $hd'$ for $h \in C(X)$. The function $fhd' = [f/(1 + f^2)]hd$ belongs to $C(\operatorname{coz} d)$ and, clearly, it extends continuously to all of $X$. Thus, $f \in \operatorname{Hom}(d')$. ∎

Thus $\mathrm{Q}(X)$ is the ring of all continuous functions on the dense open sets in $X$, and $\mathrm{Q}_{\mathrm{cl}}(X)$ is the ring of all continuous functions on the dense cozero-sets in $X$.

COROLLARY. $\mathrm{Q}(X)$ *is regular.*

PROOF. $C(X)$ is semi-prime. Apply Theorem 1.11. ∎

REMARK. It is interesting to prove this directly. If $f \in \mathrm{Q}(X)$, then $f \in C(V)$ for some dense open $V$. Define $g(x) = 1/f(x)$ for $x \in \operatorname{coz} f$, and $g(x) = 0$ for $x \in \operatorname{int} Z(f)$. Then $g$ is continuous on a dense open set, and $f^2 g = f$.

**2.7. ELABORATIONS.** If $B$ is any ring of quotients of $C(X)$, so that $B \subset \mathrm{Q}(X)$, then $B$ is a subring of a direct limit of rings of continuous functions. More explicitly, $B = \varinjlim_D (B \cap \operatorname{Hom} D)$ and $\operatorname{Hom} D \subset C(\operatorname{coz} D)$, so that $B$ is exhibited as a direct limit itself of rings of continuous functions.

We can also consider any subring $A$ of $C(X)$. Putting

$$\mathcal{V} = \bigcup_{D \in \mathcal{D}_0(A)} \{V \in \mathcal{V}_0 : V \supset \operatorname{coz} D\},$$

we have

$$\mathcal{Q}(A) = \varinjlim_{D \in \mathcal{D}_0(A)} \operatorname{Hom} D \subset \varinjlim_{D \in \mathcal{D}_0(A)} C(\operatorname{coz} D)$$
$$\xleftarrow{\tau} \varinjlim_{V \in \mathcal{V}} C(V) \subset \varinjlim_{V \in \mathcal{V}_0} C(V) = \mathrm{Q}(X),$$

where $\tau$ is the restriction homomorphism. Thus, $\mathcal{Q}(A)$ is exhibited as a subring of a homomorphic image of a subring of $\mathrm{Q}(X)$. The condition that $\tau$ be one-one is that $\operatorname{coz} D$ be dense in $V$ in each case. Hence we have:

COROLLARY.    Let $A \subset C(X)$; then $\mathcal{Q}(A) \subset Q(X)$ *if and only if all dense ideals in $A$ have dense cozero-sets.*

Finally, we see as above that every ring of quotients of $A$ is a ring of continuous functions on open sets in $X$—and, when dense ideals have dense cozero-sets, on dense open sets in $X$.

2.8   THEOREM    $Q(\mathbf{R}^n) \neq Q(\mathbf{R})$ *for $n > 1$.*

REMARK.    This result, due to J. G. Fortin and F. Rothberger, appears in Fortin's McGill thesis [Fortin (1963)]. The exposition below has benefited from some suggestions made by A. W. Hager.

We shall use the following notation: if $U$ is an open set in a space $X$, then $e_U$ denotes the idempotent in $Q(X)$ for which $\coz e_U = U$ (that is, $e_U$ is equal to 1 everywhere on $U$ and to 0 on $X - \cl U$). Also, in what follows, we use the same letter to denote a member of $Q(X)$ and any one of its representative functions defined on a dense open set in $X$.

LEMMA 1.   *Let $f \in Q(X)$ and let $e$ be an idempotent in $Q(X)$. Then the set $\cl f[\coz e]$ consists precisely of those $r \in \mathbf{R}$ such that for every $\varepsilon > 0$,*

(1) *there exists an idempotent $e_U$ in $Q(X)$ satisfying:*

  (a) $e_U \neq 0$,

  (b) $e \cdot e_U = e_U$,

  (c) $e_U \cdot [\varepsilon^2 - (f - r)^2] \geq 0$.

PROOF.    By definition, $\cl f[\coz e]$ is the set of all $r \in \mathbf{R}$ such that for every $\varepsilon > 0$,

(2) *there exists $x_0$ such that $|f(x_0) - r| \leq \varepsilon$ and $e(x_0) = 1$.* Since $(c)$ states that $|f(x) - r| \leq \varepsilon$ whenever $e_U(x) = 1$, it is clear that (1) implies (2).

Conversely, consider the set

$$U = \coz e \cap f^{\leftarrow}[(r - \varepsilon, r + \varepsilon)].$$

The corresponding idempotent $e_U$ clearly satisfies conditions $(b)$ and $(c)$. Because of (2), it satisfies $(a)$ as well. Thus, (1) holds. ∎

LEMMA 2. *If $f \mapsto f'$ is an isomorphism from $Q(X)$ onto $Q(X')$, then $\operatorname{cl} f[\operatorname{coz} e] = \operatorname{cl} f'[\operatorname{coz} e']$.*

PROOF. Any isomorphism carries **R** (the constant functions) onto **R**, identically; the proof is like that in [Gillman & Jerison (1960), 1.9 and 1I]. Isomorphisms also preserve order (since nonnegative functions are squares). The characterization of $\operatorname{cl} f[\operatorname{coz} e]$ given in Lemma 1 is thus an algebraic one, and shows that the set is preserved under isomorphisms. ∎

PROOF OF THEOREM. Assume that there exists an isomorphism $f \mapsto f'$ from $Q(\mathbf{R}^n)$ onto $Q(\mathbf{R})$, where $n > 1$. Let $i$ denote the identity function on **R** ($i(x) = x$) and let $h$ denote the preimage of $i$ in $Q(\mathbf{R}^n)$ ($h' = i$).

First we show that $h$ cannot be constant on any nonvoid open set $G$. Suppose, on the contrary, that $h[G] = \{r_0\}$. Then $\operatorname{cl} h[\operatorname{coz} e_G] = \{r_0\}$. By Lemma 2, $\operatorname{cl} i[\operatorname{coz} e'_G] = \{r_0\}$; then $\operatorname{coz} e'_G$ is a single point, which is absurd.

Accordingly, we can find a connected open set $S$ in the domain of $h$ and points $p$ and $q$ in $S$ such that $h(p) < h(q)$. Choose disjoint connected open sets $U$ and $V$ in $S$ whose closures contain both $p$ and $q$. (Here we need $n > 1$.) The connected sets $\operatorname{cl} h[U]$ and $\operatorname{cl} h[V]$ contain both $h(p)$ and $h(q)$, and hence

(3) $\qquad \operatorname{cl} h[U] \cap \operatorname{cl} h[V]$ contains the entire interval $[h(p), h(q)]$.

On the other hand, since $U \cap V = \emptyset$, we have $e_U \cdot e_V = 0$, whence $e'_U \cdot e'_V = 0$, and therefore $\operatorname{coz} e'_U \cap \operatorname{coz} e'_V = \emptyset$, from which it follows that

(4) $\qquad \operatorname{cl} \operatorname{coz} e'_U \cap \operatorname{cl} \operatorname{coz} e'_V$ is nowhere dense.

Since $\operatorname{cl} h[U] = \operatorname{cl} h[\operatorname{coz} e_U] = \operatorname{cl} \operatorname{coz} e'_U$ and, similarly,
$\operatorname{cl} h[V] = \operatorname{cl} \operatorname{coz} e'_V$, the results (3) and (4) are contradictory. This completes the proof of the theorem.

# 3. Equalities Among Various Rings of Quotients of $C(X)$

**3.1.** PRELIMINARIES. We let $Q^*(X)$ and $Q^*_{cl}(X)$ denote the subrings of all bounded functions in $Q(X)$ and $Q_{cl}(X)$, resp. Evidently, $C^*(X) \subset Q^*_{cl}(X) \subset Q^*(X)$.

*Continue* means to extend continuously; *continuation*, continuous extension.

A subspace $S$ of $X$ is said to be *C-embedded* in $X$ if every function in $C(S)$ can be continued to $X$; similarly, $S$ is $C^*$-*embedded* if every function in $C^*(S)$ can be continued to $X$. $C$-embedding implies $C^*$-embedding, but not conversely. (See [Gillman & Jerison (1960)] for a full discussion of this and related matters.)

As is well known, each (completely regular) space $X$ has a compactification $\beta X$, unique up to homeomorphism, in which $X$ is $C^*$-embedded. ($\beta X$ is known as the "Stone-Čech" compactification of $X$.) More generally, any continuous mapping from $X$ to a compact space $Y$ has a continuation from $\beta X$ into $Y$.

Since $C^*(\beta X) = C^*(X)$, we have $Q(\beta X) = Q(X)$ and $Q_{cl}(\beta X) = Q_{cl}(X)$.

The condition $C(X) = C^*(X)$ defines $X$ as *pseudocompact*. Compact spaces are pseudocompact, but not conversely. (A compact space is a Hausdorff space in which every open cover has a finite subcover.)

The device of continuing a function from an open set $U$ to a dense open set by assigning it the value 0 on $X - \operatorname{cl} U$ has several helpful consequences: if all dense open sets are $C$-embedded, or $C^*$-embedded, or pseudocompact, then *all* open sets are $C$-embedded, or $C^*$-embedded, or pseudocompact, resp.

**N** denotes the countably infinite discrete space.

**3.2.** THE SPACE $\Delta^*$. Pick an uncountable discrete space $\Delta$ and let $\Delta^* = \Delta \cup \{\infty\}$ denote its one-point compactification. This space will be useful in several examples.

Neighborhoods of $\infty$ are the complements of finite subsets of $\Delta$. It follows that all infinite sets have $\infty$ as a limit point.

Another consequence is that every nonvoid zero-set, being a $G_\delta$, meets $\Delta$. Therefore, no proper cozero-set in $\Delta^*$ is dense. In particular, $\Delta$ is not a cozero-set in $\Delta^*$.

Let $E = (e_n)$ be a sequence of distinct points in $\Delta$. If $f(e_n) = 1/n$, and $f(x) = 0$ for $x \in \Delta^* - E$, then $f$ is continuous; hence $E$ is a cozero-set in $\Delta^*$. A function on $E$ that assumes distinct constant values on complementary infinite sets cannot be continued to the point $\infty$; consequently, the cozero-set $E$ is not $C^*$-embedded in $\Delta^*$. Likewise, $\Delta$ is not $C^*$-embedded.

**3.3.   $Q(X)$ AND $Q_{cl}(X)$.** The equality $Q(X) = Q_{cl}(X)$ means that every continuous function on a dense open set in $X$ can be defined on a dense cozero-set (2.6); similarly, $Q^*(X) = Q^*_{cl}(X)$ means that every such bounded function can be so defined.

THEOREM.    $Q(X) = Q_{cl}(X)$ *if and only if* $Q^*(X) = Q^*_{cl}(X)$.

PROOF. The necessity is trivial. Conversely, given $h \in C(V)$, where $V$ is dense and open in $X$, put $h' = 1/(1+h^2)$; then $h' \in C^*(V)$ and $hh' \in C^*(V)$. By hypothesis, there exist $f, g \in C(X)$ such that $h$ is defined on $\coz f$ and $hh'$ on $\coz g$. Since $h'$ has no zeros, $h = (hh')/h'$ is defined on $\coz f \cap \coz g$. ∎

If $X$ is a metric space, then $Q(X) = Q_{cl}(C)$—for in a metric space, every open set is a cozero-set.

The space $\Delta^*$ satisfies $Q \neq Q_{cl}$, since $\Delta$ is not $C^*$-embedded in $\Delta^*$, the only dense cozero-set in $\Delta^*$.

**3.4.   $Q_{cl}(X)$ AND $C(X)$.**   The equality $Q_{cl}(X) = C(X)$ means that every dense cozero-set in $X$ is $C$-embedded, hence closed (since $1/f \in C(\coz f)$.) Therefore $Q_{cl}(X) = C(X)$ if and only if no proper cozero-set in $X$ is dense.

$Q^*_{cl}(X) = C^*(X)$ means that every dense cozero-set in $X$ is $C^*$-embedded.

Obviously, $Q_{cl} = C$ implies $Q^*_{cl} = C^*$. The space $\beta\mathbf{N}$ satisfies $Q^*_{cl} = C^*$ but $Q_{cl} \neq C$ (since $\mathbf{N}$ is a dense cozero-set).

$Q_{cl}(X) = C(X)$ does not imply that *all* cozero-sets in $X$ are $C^*$-embedded (i.e., that $X$ is an "$F$-space" [Gillman & Jerison (1960), 14.25]); the space $\Delta^*$ is an example (supplied by Henriksen, simpler than the original one).

**3.5.   $Q(X)$ AND $C(X)$.**   The equality $Q(X) = C(X)$ means that every dense open set is $C$-embedded, whence *every* open set is $C$-embedded. This requires that $X$ be *extremally disconnected* (every open set is $C^*$-embedded [Gillman & Jerison (1960), 1H]) and a *P-space* (every cozero-set is $C$-embedded [Gillman & Jerison (1960), 14.29]). Since $C^*$-embedding in a $P$-space implies $C$-embedding (as follows from [Gillman & Jerison (1960), 1.18]), the converse holds as well.

It is known that extremally disconnected $P$-spaces of reasonable cardinal (not large enough to admit a 2-valued measure) are discrete [Gillman & Jerison (1960), 12H]; for "practical" purposes, then, $Q(X) = C(X)$ if and only if $X$ is discrete.

$Q^*(X) = C^*(X)$ means that every dense open set is $C^*$-embedded, i.e., that $X$ is extremally disconnected.

The space $\beta\mathbf{N}$ satisfies $Q^* = C^*$ but $Q \neq C$.

**3.6.  Q(X) AND Q*(X).** The equality $Q(X) = Q^*(X)$ means that every dense open set is pseudocompact, whence *every* open set is pseudocompact. This requires that $X$ be finite; for else some $f \in C(X)$ assumes arbitrarily small positive values [Gillman & Jerison (1960), 3L], whence $a/f \in C(\operatorname{coz} f)$ is unbounded.

$Q_{cl}(X) = Q_{cl}^*(X)$ means that every dense cozero-set is pseudocompact. Hence $X$ is pseudocompact and no proper cozero-set is dense (since $1/f \in C(\operatorname{coz} f)$ is bounded away from 0, whence $\operatorname{coz} f$ is closed).

$\Delta^*$ is an infinite compact space containing no proper dense cozero-set (3.2). $\beta\mathbf{N} - \mathbf{N}$ is another such space [Gillman & Jerison (1960)].

**3.7.  DENSE SUBSPACES.**  Let $S$ be a dense subspace of $Y$. The condition that $Q(S) = Q(Y)$ is that every continuous function on a dense open set in $S$ can be defined on some dense open set in $Y$. Similarly, $Q_{cl}(S) = Q_{cl}(Y)$ means that every continuous function on a dense cozero-set in $S$ can be defined on a dense cozero-set in $Y$.

In 3.12 we shall present an example in which $Q(S) \neq Q(Y)$ as well as an example where $S$ is not open but $Q(S) = Q(Y)$.

$Q(S) = Q(Y)$ does not imply $Q_{cl}(S) = Q_{cl}(Y)$, even for $S$ open. For instance, $Q_{cl}(\Delta) = Q(\Delta) = Q(\Delta^*) \neq Q_{cl}(\Delta^*)$.

The above conditions state that for each dense open (or cozero-) set $U$ in $S$ and for each $h \in C(U)$, there exist a dense open (or cozero-) set $W$ in $Y$ and a function $g \in C(W)$ such that $g|_{W \cap U} = h|_{W \cap U}$. It is not required that $W \supset U$; indeed, we shall present an example in 3.12 in which for a particular $h$, no $W \supset U$ exists.

First, however, we consider the important case $S = X$, $Y = \beta X$. We remark that in general, $X$ is not open in $\beta X$; in fact, $X$ is open in $\beta X$ if and only if $X$ is locally compact [Gillman & Jerison (1960)].

**3.8  THEOREM**  *Every continuous function on a dense open set $V$ in $X$ can be continued to an open set $(\supset V)$ in $\beta X$. Every continuous function on a dense cozero-set in $X$ can be continued to a cozero-set in $\beta X$.*

PROOF. Given $h \in C(V)$, define $D = h^{-1}C(\beta X)$; then $h \in \operatorname{Hom}(D, C(\beta X))$. By Lemma 2.5, $h \in C(\operatorname{coz}_{\beta X} D)$. By (4) of 2.3, $\operatorname{coz} D \supset V$. In case $V$ is a cozero-set $\operatorname{coz}_X f$, we define $d = f/(1 + f^2 + h^2)$; then $d \in D$, $V \subset \operatorname{coz}_{\beta X} d \subset \operatorname{coz} D$, and $h \in C(\operatorname{coz} d)$.

COROLLARY.    $C(X) = \varinjlim_W C(W)$, $W$ ranging over all cozero-sets in $\beta X$ that contain $X$.

Explicitly, each $h \in C(X)$ can be continued to $\operatorname{coz} d$, where $d = 1/(1 + h^2)$. In the notation of [Gillman & Jerison (1960), 8B], $\operatorname{coz} d$ is the space $\upsilon_h X$; thus,

$$C(X) = \varinjlim_{h \in C(X)} C(\upsilon_h X).$$

**3.9.** IRREDUCIBLE FUNCTIONS.     A continuous function defined on a subspace of $X$ will be called *irreducible* if it cannot be continued further in $X$.

Every continuous function on a dense set $S$ has an irreducible continuation; the union $h$ of all its continuations. (To see that $h$ is continuous, note that each $x \in \operatorname{dom} h$ has a neighborhood (in $X$) within which $h(x)$ is near $h(s)$ for $s \in S$ and hence near $h(y)$ for $y \in \operatorname{dom} h$.)

Every continuous function on an open set $U$ has a continuation to a dense open set (obtained by assigning it the value 0 throughout $X - \operatorname{cl} U$) and hence has an irreducible continuation.

A continuous function $f$ can fail to have an irreducible continuation. For this, it is sufficient that $f$ be defined on a closed set $E$ whose complement $D$ is discrete and that $f$ have no continuation to all of $X$: for then if $f'$ is any continuation of $f$, there is an isolated point outside its domain, and $f'$ can be continued further to that point. As an example, let $E$ be the real axis, $D$ the rational points in the upper half plane, and $X = E \cup D$, with the relative topology of the plane enlarged as follows: each point of $D$ is isolated, and the neighborhoods of $a = (x, 0)$ include all sets $\{b \in D : |b - (x, r)| < r\} \cup \{a\}$ for $r > 0$. One can check that $X$ is completely regular (cf. [Gillman & Jerison (1960), 3K]). Clearly, the countable set $D$ is dense in $X$, and $E$ is discrete; hence $X$ admits just $c$ continuous functions, while $E$ admits $2^c$, so that some function on $E$ has no continuation to $X$. (For alternative examples, see [Gillman & Jerison (1960), 5I or 6Q].)

3.10   THEOREM

(1) *The domain of an irreducible function is a dense $G_\delta$.*

(2) *Conversely, if $X$ is locally connected, then each dense, countable intersection of cozero-sets is the domain of some bounded irreducible function.*

PROOF. (1). This result is familiar from classical analysis. Let $f \in C(S)$ be irreducible; then $S$ is dense. To prove that $S$ is a $G_\delta$, one first defines $\operatorname{osc} f$ on $X$:

$$\operatorname{osc} f(x) = \inf_U [\sup f(U \cap S) - \inf f(U \cap S)],$$

$U$ ranging over all neighborhoods of $x$. Then $\operatorname{osc} f$ is upper semi-continuous and

$$S = \{x \in X : \operatorname{osc} f(x) = 0\} = \bigcap_{n>0} \{x \in X : \operatorname{osc} f(x) < 1/n\},$$

so that $S$ is a $G_\delta$.

(2). Let $S = \bigcap_{n>0} \operatorname{coz}_X g_n$ be dense. Define $f_n = \sin(1/g_n)$ on $\operatorname{coz} g_n$; then $f_n \in C(\operatorname{coz} g_n)$ and $|f_n| \leq 1$. Next, we show that for each $x \in Z(g_n)$, $\operatorname{osc} f_n(x) = 2$. Let $U$ be any connected neighborhood of $x$. Then $g_n(U)$ is a connected set in $\mathbf{R}$, containing 0 (because $x \in U$), and containing points $\neq 0$ (because the dense set $\operatorname{coz} g_n$ meets $U$).

Therefore $g_n(U \cap \text{coz } g_n)$ contains an open interval having 0 as an endpoint. This implies that osc $f_n(x) = 2$.

Now define $f = \sum_{n>0} 3^{-n} f_n$ on $S$; then $f \in C^*(S)$. To prove that $f$ is irreducible, we show that osc $f(s) > 0$ for every $x \notin S$. Given such $x$, let $m$ be the smallest integer such that $x \in Z(g_m)$. Then

$$\begin{aligned} \text{osc } f(x) &\geq 3^{-m} \text{osc } f_m(x) - \sum_{n>m} 3^{-n} \text{osc } f_n(x) \\ &\geq 2 \cdot (3^{-m} - \sum_{n>m} 3^{-n}) = 3^{-m}. \end{aligned}$$

∎

**3.11 COROLLARY** *In a locally connected metric space, each dense $G_\delta$ is the domain of some bounded irreducible function.*

PROOF. Every open set in a metric space is a cozero-set. ∎

REMARK. In (1), we cannot derive the stronger conclusion that the domain is a dense, countable intersection of *cozero*-sets. For example, $\Delta$ is not $C^*$-embedded in $\Delta^*$, the only cozero-set containing it.

In (2), the hypothesis of local connectedness cannot be wantonly discarded. In the space $\beta\mathbf{N}$, for example, $\mathbf{N}$ is a dense cozero-set, but every bounded function on $\mathbf{N}$ continues to all of $\beta\mathbf{N}$.

**3.12. EXAMPLES.** (See 3.7.) (1). *$S$ dense in $Y$ but $Q(S) \neq Q(Y)$.* Let $Y = \mathbf{R}$ and let $S$ be the subspace of irrationals. Since $S$ is a dense $G_\delta$ in the locally connected metric space $\mathbf{R}$, it is the domain of some irreducible function $h$; then $h$ is continuous on the dense open set $S$ in $S$. If $V$ is dense and open in $\mathbf{R}$, then $h|_{V \cap S}$ cannot be continued to any rational point—else $h$ itself could be continued to that point (see 3.9). Since $V$ contains rationals, no function in $C(V)$ can agree with $h$ on $V \cap S$.

(2). *$S$ dense (but not open) in $Y$, with $Q(S) = Q(Y) = Q_{cl}(S) = Q_{cl}(Y)$, but for some dense cozero-set $U$ in $S$, some $h \in C(U)$ has no continuation to a dense open set $(\supset U)$ in $Y$.* Let $Y$ be a closed disc in the plane and $S$ the set int $Y$ plus a boundary point of $Y$. Then $S$ is dense in the metric space $Y$. To show that $Q(S) = Q(Y)[= Q_{cl}(S) = Q_{cl}(Y)]$, let $V$ be dense and open in $S$ and let $f \in C(V)$. Then $V \cap \text{int } Y$ is dense and open in $Y$ and $f \in C(V \cap \text{int } Y)$.

Now observe that $S$ is a dense $G_\delta$ in the locally connected metric space $Y$. Let $h$ be an irreducible function with domain $S$; then $h$ is continuous on the dense open set $U = S$ in $S$. But $h$ has no continuation to an open set in $Y$ (note that $S$ itself is not open in $Y$), because it has no proper continuation to anything in $Y$.

# 4. Various Operations on Q(X)

4.1  DEFINITION   *If $\mathcal{S}(X)$ is any filter base of dense subsets of $X$, then*
$$C[\mathcal{S}] = \varinjlim_{S \in \mathcal{S}} C(S)$$
*is defined as a ring of continuous functions (2.4). In this notation, $Q(X) = C[\mathcal{V}_0(X)]$. The subring of bounded functions is denoted $C^*[\mathcal{S}]$. Clearly,*
$$C^*[\mathcal{S}] = \varinjlim_{S \in \mathcal{S}} C^*(S).$$
*As in 2.3, we see that $C[\mathcal{S}]$ is a ring of quotients of $C^*[\mathcal{S}]$.*

*Let $f \in C^*[\mathcal{S}]$; then $f \in C^*(S)$ for some $S \in \mathcal{S}$. The sup norm $\|f\| = \sup_{s \in S} |f(s)|$ is determined by the values of $f$ on any dense subset of $S$. Hence we may define the sup norm of $f$ in $C^*[\mathcal{S}]$ to be its norm in $C(S)$ and denote it by the symbol $\|f\|$. The sup norm induces a natural metric $(f, g) \to \|f - g\|$ on $C^*[\mathcal{S}]$. Then $C^*[\mathcal{S}]$ becomes a normed algebra over $\mathbf{R}$ (not necessarily complete).*

*An equivalent metric on $C^*[\mathcal{S}]$ is*
$$(f, g) \to \|\frac{f - g}{1 + |f - g|}\|,$$
*and the same expression defines a metric on all of $C[\mathcal{S}]$. (Note, however, that multiplication in $C[\mathcal{S}]$ is not in general continuous—not even in $C(X)$.) Because the two metrics on $C^*[\mathcal{S}]$ are equivalent, we may deal with convergence in $C[\mathcal{S}]$ in terms of the natural metric or sup norm. Evidently, $C^*[\mathcal{S}]$ is a closed set in the metric space $C[\mathcal{S}]$. The metric completion of $C$ is denoted $\overline{C}$.*

*Each ring $C[\mathcal{S}]$ is also a lattice with respect to the pointwise definition of order. If every nonvoid subset with an upper bound in $C$ has a supremum in $C$, then $C$ is* Dedekind-complete. *If $C$ is Dedekind-complete, $B \subset C$, and every element of $C$ is the supremum of some subset of $B$, then $C$ is the* Dedekind completion *of $B$.*

4.2  THEOREM   *Let $\mathcal{S}$ be a filter base of dense sets in $X$. The following is a necessary and sufficient condition that $C[\mathcal{S}]$ be rationally complete:*

*Let $(E_i)_{i \in I}$ be any family of disjoint open subsets in $X$ whose union is dense in $X$, let $(S_i)_{i \in I}$ be a subfamily of $\mathcal{S}$, and let*
$$T = \bigcup_{i \in I} (E_i \cap S_i)$$

(whence $T$ is dense in $X$); then $C(T) \subset C[\mathcal{S}]$.

PROOF. *Necessity.* Given $g \in C(T)$, Let $D = g^{-1}C[\mathcal{S}]$; then $D$ is an ideal of $C[\mathcal{S}]$ and $g \in \mathrm{Hom}\, D$. We show that $S$ is dense, from which it will follow that $g \in \mathcal{Q}(C[\mathcal{S}]) = C[\mathcal{S}]$. Suppose $0 \neq h \in C[\mathcal{S}]$. Then $\mathrm{coz}\, h$ meets the dense open set $\bigcup_{i \in I} E_i$, hence meets some $E_i$, and therefore meets $E_i \cap S_i$ (since $S_i \cap \mathrm{dom}\, h$ is dense). So $h(x) \neq 0$ for some $x \in E_i \cap S_i$. Choose $d \in C(X)$ to vanish on a neighbourhood of $X - E_i$ but not at $x$. The function $f$ equal to $gd$ on $E_i \cap S_i$ and $0$ on $X - E_i$ is continuous on $(E_i \cap S_i) \cup (X - E_i)$, and it agrees with $gd$ on the dense set $T$. But $f \in C(S_i) \subset C[\mathcal{S}]$. Therefore, $gd \in C[\mathcal{S}]$, i.e., $d \in D$. Obviously, $hd \neq 0$.

*Sufficiency.* Let $\varphi \in \mathrm{Hom}\, D$, where $D$ is a dense ideal in $C[\mathcal{S}]$; we are to show that $\varphi \in C[\mathcal{S}]$. Represent each element of $D$ by a specific function $d$ defined on some member of $\mathcal{S}$. Evidently, $\bigcup_d \mathrm{coz}\, d$ is dense in $X$ (else $D$ would have a nonzero annihilator in $C(X)$). Pick $H_d$ open in $X$ such that $\mathrm{coz}\, d = H_d \cap \mathrm{dom}\, d$; then $\bigcup_d H_d$ is open and dense in $X$. By Zorn's Lemma, there exist disjoint open sets $E_d \subset H_d$ such that $\bigcup_d E_d$ is also dense.

Now for each $d$, pick a representative $\varphi(d)$, and define $S_d = \mathrm{dom}\, d \cap \mathrm{dom}\, \varphi(d)$ and

$$T = \bigcup_d (E_d \cap S_d).$$

Then $S_d \in \mathcal{S}$ and, by hypothesis, $C(T) \subset C[\mathcal{S}]$.

Consider any $x \in T$. There is a unique $d$ such that $x \in E_d \cap S_d$; define

$$g(x) = \frac{\varphi(d)(x)}{d(x)}.$$

This defines $g$ on all of $T$. Since $g$ is continuous on each of the open subsets $E_d \cap S_d$ of $T$, it is continuous on $T$. So $g \in C[\mathcal{S}]$.

Finally, consider any $d \in D$. For every $x \in T \cap S_d$,

$$g(x)d(x) = \frac{\varphi(d')(x)}{d'(x)} d(x) = \varphi(d)(x) \qquad (x \in E_{d'} \cap S_{d'}).$$

Therefore $gd$ agrees with $\varphi(d)$ on a dense subset; hence $gd = \varphi(d)$. Consequently, $g = \varphi$. ∎

REMARK. The inclusion $C(T) \subset C[\mathcal{S}]$ does not imply $T \in \mathcal{S}$, even when $\mathcal{S}$ is a filter. Let **Q** denote the rationals[1] and let $\mathcal{S}$ be all sets in **R** that contain $G_\delta$'s containing **Q**. Every function in $C(\mathbf{Q})$ has an irreducible extension to a dense $G_\delta$ (3.9 and 3.10). But $\mathbf{Q} \notin \mathcal{S}$.

---

[1] The original used $P$ for the rationals, presumably because they lacked a bold font.

**4.3.** LOCALLY CONSTANT FUNCTIONS. We shall say that $f \in C(S)$ is *locally constant* provided that $f^{\leftarrow}(r)$ is open in $S$ for each $r \in \mathbf{R}$. The locally constant functions in $C(S)$ form a ring $L(S)$. We define

$$\mathrm{Q}_L(X) = \varinjlim\nolimits_{V \in \mathcal{V}_0(X)} L(V).$$

Then $\mathrm{Q}_L(X) \subset \mathrm{Q}(X)$. The set of functions in $C(S)$ whose range is finite is denoted by $F(S)$. We define

$$\mathrm{Q}_F(X) = \varinjlim\nolimits_{V \in \mathcal{V}_0(X)} F(V).$$

Since $F(V) \subset L^*(V)$, we have $\mathrm{Q}_F(X) \subset \mathrm{Q}_L^*(X)$. Theorem 3.8 leads to a proof that $\mathrm{Q}_L(\beta X) = \mathrm{Q}_L(X)$ and $\mathrm{Q}_F(\beta X) = \mathrm{Q}_F(X)$.

The ring $\mathrm{Q}_L(X)$ is rationally complete. To see this, notice first that if $D$ is a dense ideal in $\mathrm{Q}_L(X)$, then $\operatorname{coz} D$ is dense in $X$; for if $\operatorname{coz} D$ misses a nonvoid open set $G$, then the characteristic function $\chi_G$ on $G \cup (X-)\operatorname{cl} G$ belongs to $\mathrm{Q}_L(X)$, is not zero, and annihilates $D$. The result can now be inferred from the proof of sufficiency in Theorem 4.2; more directly, we may argue as follows. Let $\varphi \in \operatorname{Hom} D$ be given. For each $r \in \mathbf{R}$, define $g(x) = r$ for all $x$ in the open set

$$\bigcup_{d \in D} \bigcup_{0 \neq s \in \mathbf{R}} [d^{\leftarrow}(s) \cap \varphi(d)^{\leftarrow}(rs)].$$

Then $g$ is defined on a dense open set (denseness following from the denseness of $\operatorname{coz} D$); evidently, $g$ is locally constant. Finally, one verifies without difficulty that for each $d \in D$, $gd$ agrees with $\varphi(d)$ on a dense set.

$\mathrm{Q}_L(X)$ is a rational extension of $\mathrm{Q}_F(X)$ and hence is its rational completion. For consider any $h \in \mathrm{Q}_L(X)$, $h \neq 0$. Then $0 \neq h \in L(V)$ for some dense open $V$, and there exists $r \neq 0$ for which $G = h^{\leftarrow}(r)$ is nonvoid. The function $\chi_G$ on $V$ belongs to $F(V)$ and $0 \neq h \cdot \chi_G \in F(V) \subset \mathrm{Q}_F(X)$.

4.4 LEMMA   *If $\mathcal{S} \supset \mathcal{V}_0(X)$, then $\mathrm{Q}_L(X)$ and hence $\mathrm{Q}(X)$ are metrically dense in $C[\mathcal{S}]$, and $\mathrm{Q}_F(X)$ and $\mathrm{Q}^*(X)$ in $C^*[\mathcal{S}]$. Hence $C[\mathcal{S}] = \mathrm{Q}_L(X) + C^*[\mathcal{S}] = \mathrm{Q}(X) + C^*[\mathcal{S}]$.*

PROOF. The first assertion implies the last. It also implies that $\mathrm{Q}_L^*(X)$ is dense in $C^*[\mathcal{S}]$; since, obviously, $\mathrm{Q}_F(X)$ is dense in $\mathrm{Q}_L^*(X)$, the second assertion follows.

To prove the first, let $S \in \mathcal{S}$, $g \in C(S)$, and an integer $n > 0$ be given; we are to find $V \in \mathcal{V}_0$ and $f \in L(V)$ such that $\|g - f\| \leq 1/n$. Let $E$ denote the set of all integer multiples of $1/n$. For $e \in E$, the set

$$U_e = g^{\leftarrow}(e, e + 1/n) \cup \operatorname{int}_S g^{\leftarrow}(e)$$

is open in $S$; choose $V_e$ open in $X$ such that $U_e = V_e \cap S$. Since $E$ is discrete, $\operatorname{int} g^{\leftarrow}(E) = \bigcup_{e \in E} \operatorname{int} g^{\leftarrow}(e)$.

Next, $U = \bigcup_{e \in E}$ is dense in $S$; for if $s \in S - U$, then $s \notin \operatorname{int} g^{\leftarrow}(E)$, i.e., every neighbourhood of $s$ meets $g^{\leftarrow}(\mathbf{R} - E) \subset U$.

Since $S$ is dense in $X$, so, then, is $U$ and hence so is $V = \bigcup_{e \in E} V_e$.

Next, the open sets $V_e$ are disjoint; for if $e \neq e'$, then $V_e \cap V_{e'} \cap S = U_e \cap U'_e = \emptyset$, whence $V_e \cap V_{e'} = \emptyset$, since $S$ is dense. Therefore,

$$f(v) = e \qquad (v \in V_e, e \in E)$$

defines $f$ as a locally constant function on $V \in \mathcal{V}_0$. Clearly, $\|g - f\| \leq 1/n$. ∎

**4.5 LEMMA** *If $\mathcal{S}$ is closed under countable intersection, then $C[\mathcal{S}]$ is metrically complete.*

PROOF. Let $(g_n)$ be a Cauchy sequence in $C[\mathcal{S}]$. Each $g_n$ belongs to $C(S)$, where $S = \bigcap_n \operatorname{dom} g_n \in \mathcal{S}$. Since $C(S)$ is complete, $\lim_n g_n$ exists in $C(S)$. ∎

A *Baire* space is one in wich countable intersections of dense open sets are dense. Locally compact spaces and complete metric spaces are Baire spaces, as is well known.

$\mathcal{G}_0(X)$ denotes the family of all dense $G_\delta$'s in $X$. Clearly, in a Baire space, $\mathcal{G}_0$ is closed under countable intersection.

**4.6 THEOREM** $\overline{Q}(X) = C[\mathcal{G}_0(\beta X)]$; *and if $X$ is a Baire space, then $\overline{Q}(X) = C[\mathcal{G}_0(X)]$.*

PROOF. Recall that $Q(X) = Q(\beta X) = C[\mathcal{V}_0(\beta X)]$. Since $\beta X$ is compact, $\mathcal{G}_0(\beta X)$ is closed under countable intersection. By Lemma 4.4, $Q(\beta X)$ is dense in $C[\mathcal{G}_0(\beta X)]$; and by Lemma 4.5, the latter is metrically complete. This establishes the first assertion; the second follows similarly. ∎

**4.7. THE RINGS $\overline{Q}(X)$ AND $\overline{Q}^*(X)$.** We have just proved that $\overline{Q}(X)$ is a ring of continuous functions. We are therefore invited to consider its subring $\overline{Q}^*(X)$ of bounded functions. From the relation $Q^*(X) \subset \overline{Q}^*(X)$, we get $Q^*(X) = Q(X) \cap \overline{Q}^*(X)$.

Since $\overline{Q}^*(X)$ is a closed set in $\overline{Q}(X)$, it is complete; therefore $\overline{Q}^*(X)$ is a Banach algebra over **R**. Evidently, $Q^*(X)$ is dense in $\overline{Q}^*(X)$; therefore,

$$\overline{Q^*}(X) = \overline{Q}^*(X).$$

Since $Q_L(X) = Q_L(\beta X)$, $Q_L(X)$ is dense in $\overline{Q}(X)$, $Q_F(X)$ is dense in $\overline{Q}^*(X)$, and $\overline{Q}(X) = Q_L + \overline{Q}^*(X) = Q(X) + \overline{Q}^*(X)$ (Lemma 4.4).

In general, $\overline{Q}(X) \neq Q(X)$. For example, $\overline{Q}(\mathbf{R}) \neq Q(\mathbf{R})$, as follows from Example (1) in 3.12. Since $C(\mathbf{R})$ is metrically complete, this can be written $\overline{\mathcal{Q}}(C(\mathbf{R})) \neq \mathcal{Q}(\overline{C}(\mathbf{R}))$; thus the operations of metric completion and rational completion need not commute. When acting on $Q^*(X)$, however, they do commute, as we now show.

4.8 THEOREM $\overline{Q}(X)$ *is rationally complete. Hence* $\overline{Q}(X) = \mathcal{Q}(\overline{Q}^*(X))$ *and* $\overline{Q}^*(X) = \mathcal{Q}^*(\overline{Q}^*(X))$.

PROOF. We verify the condition of Theorem 4.2—in fact, the following stronger result: if $E_i$ are open and disjoint and $S_i$ are $G_\delta$'s, then $T = \bigcup_i (E_i \cap S_i)$ is a $G_\delta$. We have $S_i = \bigcap_n V_{in}$, $V_{in}$ open ($n = 1, 2, \cdots$). Define

$$T_n = \bigcup_i (E_i \cap V_{i\,1} \cap \cdots \cap V_{i\,n});$$

then $T_n$ is open. Because the $E_i$ are disjoint, $T = \bigcap_n T_n$. So $T$ is a $G_\delta$. ∎

4.9 THEOREM $\overline{Q}(X)$ *is the Dedekind completion of* $Q_L(X)$ *and hence of* $Q(X)$. *Likewise,* $\overline{Q}^*(X))$ *is the Dedekind completion of* $Q_F(X)$ *and of* $Q^*(X)$.

PROOF. We observe that metric denseness implies order-denseness: if $g_n \to h$, then $r_n = \|h - g_n\| \to 0$ whence $h = \sup_n(g_n - r_n)$; therefore every element of $\overline{Q}$ is the supremum of some subset of $Q_L$ and every element of $\overline{Q}^*$ is the supremum of some subset of $Q_F$ (4.7). Next, if $\overline{Q}$ is Dedekind-complete, then so is $\overline{Q}^*$. Consequently, the last half of the theorem is implied by the first.

Now let $P$ be any nonvoid subset of $\overline{Q} = C[\mathcal{G}_0(\beta X)]$ with an upper bound $u$. In seeking a supremum of $P$, we may assume that $P \subset Q_L$ and $u \in Q_L$; furthermore, we may assume that each $f \in P$ takes irrational values only, while $u$ takes rational values only.

For integers $n \geq 0$ and $k \geq 0$, and for rational $q$, define

$$V^n_{q,k} = u^{\leftarrow}(q) \cap \bigcup_{f \in P} f^{\leftarrow}(q - \frac{k+1}{2^n}, q - \frac{k}{2^n}).$$

Then $V^n_{q,k}$ is open in $\beta X$ and $\bigcup_k V^n_{q,k} = u^{\leftarrow}(q) \cap \bigcup_f \mathrm{dom}\, f$.

Next, define

$$W^n_{q,k} = V^n_{q,k} - cl \bigcup_{1 < k} W^n_{q,i} \qquad [= V^n_{q,k} - cl \bigcup_{i < k} V^n_{q,i}].$$

Then $\bigcup_k W^n_{q,k}$ is dense in $\bigcup_k V^n_{q,k}$. Hence the set

$$W^n = \bigcup_q \bigcup_k W^n_{q,k}$$

is dense in $\bigcup_q \bigcup_k V^n_{q,k} = \mathrm{dom}\, u \cap \bigcup_f \mathrm{dom}\, f$. Consequently, $W^n$ is open and dense in $\beta X$.

Now put $S = \bigcap_n W^n$. Then $S$ is a dense $G_\delta$ in $\beta X$: $S \in \mathcal{G}_0(\beta X)$.

For each $n$ and for each $s \in S$, there are unique $q$ and $k$ such that $s \in W^n_{q,k}$; define

$$g_n(s) = q - \frac{k}{2^n} \qquad\qquad (s \in W^n_{q,k}).$$

This defines $g_n$ on $S$. Evidently, $g_n \in L(S) \subset C(S)$.

Now let $n \geq 0$ be given and consider any $s \in S$. There exist unique $q$ and $k$ for which $s \in W_{q,k}^n$. Since
$$W_{q,k}^n \subset V_{q,k}^n - \bigcup_{i<k} V_{q,i}^n,$$
we have
$$g_n(s) = q - \frac{k}{2^n} > f(s) \qquad \text{for all } f \in P$$
and
$$q - \frac{k+1}{2^n} < f_0(s) \qquad \text{for at least one } f_0 \in P.$$
Therefore, letting $g = \text{pointwise sup } P$, we have
$$0 \leq g_n(s) - g(s) \leq \frac{1}{2^n}.$$
Hence, $(g_n)$ converges uniformly to $g$, $g \in C(S)$, and $g = \sup P$. ∎

**4.10 LEMMA** *Given $S \subset X$ and $h \in C(S)$, with $h \geq 0$, there exists $B \subset C^*(X)$ such that $h(s) = \sup_{f \in B} f(s)$ for all $s \in S$.*

PROOF. For $s \in S$ and $\epsilon > 0$, let $V_{s,\epsilon}$ be a neighbourhood of $s$ in $X$ such that
$$\sup h(V_{s,\epsilon} \cap S) - \inf h(V_{s,\epsilon} \cap S) < \epsilon.$$
Choose $g_{s,\epsilon} \in C^*(X)$ such that $g_{s,\epsilon}(s) = 1$, $g(X - V_{s,\epsilon}) = 0$ and $0 \leq g_{s,\epsilon} \leq 1$; and define
$$f_{s,\epsilon}(x) = \inf h(V_{s,\epsilon} \cap S) \cdot g_{s,\epsilon}(x) \qquad\qquad (x \in X).$$
Then $f_{s,\epsilon} \in C^*(X)$, $f_{s,\epsilon}|S \leq h$, and
$$f_{s,\epsilon}(s) = \inf h(V_{s,\epsilon} \cap S) \geq h(s) - \epsilon.$$
Hence for all $s \in S$,
$$h(s) - \epsilon \leq f_{s,\epsilon}(s) \leq \sup_{t \in S} f_{t,\epsilon}(s) \leq h(s).$$
Therefore, $h(s) = \sup_{t,\epsilon} f_{t,\epsilon}(s)$. ∎

**4.11 THEOREM** *The Dedekind completion of $C(X)$ is the ring of all "C-bounded" functions in $\overline{Q}(X)$—i.e., all $h \in \overline{Q}(X)$ for which there exists $g \in C(X)$ with $|h| \leq g$. Likewise the Dedekind completion of $C^*(X)$ is $\overline{Q}^*(X)$.*

PROOF. Clearly, the ring of $C$-bounded functions is Dedekind-complete, since $\overline{Q}$ is. If $h \in \overline{Q}$ is $C$-bounded (resp. bounded), then $h + g \geq 0$ for some $g \in C$ (resp. $C^*$). By the lemma, there exists $B \subset C^*$ such that $h + g = \sup B$. Then $h = \sup_{f \in B}(f - g)$. ∎

**4.12 Corollary** [Stone-Nakano] *$X$ is extremally disconnected if and only if $C(X)$ is Dedekind complete, and if and only if $C^*(X)$ is Dedekind-complete.*

Proof. The result is well-known (see, e.g., [Gillman & Jerison (1960), 3N and 6M]), but here we have another proof. Since $C^* \subset Q^* \subset \overline{Q}^*$ and $C^*$ is metrically complete, $C^* = Q^*$ if and only if $C^* = \overline{Q}^*$—that is $X$ is extremally disconnected (3.5) if and only if $C^*$ is Dedekind-complete (4.11). Next, Dedekind completeness of $C$ obviously implies that of $C^*$. Conversely, assume $\overline{Q}^* = C^*$, and let $h \in \overline{Q}$ be $C$-bounded. Then $|h| + 1 \leq g$ for some $g \in C$, whence $f = h/g \in \overline{Q}^* = C^*$ so that $h = fg \in C$. ∎

# 5. Maximal ideal spaces of $C[\mathcal{S}]$ and $C^*[\mathcal{S}]$

**5.1.** LATTICE HOMOMORPHISMS. Let $A$ be a ring, $I$ an ideal in $A$. The image of $a$ ($\in A$) in $A/I$ under the canonical homomorphism is denoted by $I(a)$. When $A$ and $A/I$ are lattices and the canonical homomorphism preserves the lattice operations, then it is called a lattice homomorphism. In particular, under a lattice homomorphism, $|I(a)| = I(|a|)$. (By definition, $|a| = a \vee -a$. For more details, see 9.6.)

$\mathfrak{M}(A)$ denotes the set of all maximal ideals in $A$.

Let $\mathcal{S}$ be a filter base of dense sets in a space $X$. For brevity, we shall write

$$C = C[\mathcal{S}],\ C^* = C^*[\mathcal{S}]$$

$$\mathcal{M} = \mathfrak{M}(C),\ \mathcal{M}^* = \mathfrak{M}(C^*).$$

LEMMA. *All homomorphisms $C \to C/M$ ($M \in \mathcal{M}$) and $C^* \to C^*/M$ ($M \in \mathcal{M}^*$) are lattice homomorphisms, and the residue class fields are totally ordered fields containing $\mathbf{R}$ as the image of the constant functions. Moreover, $C^*/M = \mathbf{R}$.*

PROOF. As in [Gillman & Jerison (1960), 5.5 and 5.8].

In general, the fields $C/M$ are non-archimedean, i.e., they contain infinitely large elements ($> n$ for every integer $n$) and infinitely small elements ($> 0$ but $< \frac{1}{n}$ for every positive integer $n$).

**5.2 THEOREM** *For $f \in C^*$, $\|f\| = \sup_{M \in \mathcal{M}^*}|M(f)|$, so that $\|f\|$ is determined algebraically. Consequently, any isomorphism of $C^*$ to a ring of functions is an isometry.*

PROOF. Since $|f| \leq \|f\|$, each $|M(f)| \leq \|f\|$ and therefore, $s = \sup_M|M(f)| \leq \|f\|$. On the other hand, for any $r \leq \|f\|$, there exists $x$ with $|f(x)| > r$. Then $(|f| - r) \wedge 0$ vanishes on a neighborhood of $x$; hence it is a non-unit of $C^*$ and therefore belongs to some maximal ideal $M$; then $|M(f)| \geq r$, whence $s \geq r$. This implies that $s \geq \|f\|$. Consequently, $s = \|f\|$.

**5.3.   THE STONE TOPOLOGY.** The *Stone* topology on $\mathfrak{M}(A)$ has for basic open sets those sets of the form
$$\Gamma(a) = \{M \in \mathfrak{M}(A) : M(a) \neq 0\} \qquad (a \in A).$$

On $\mathcal{M}$ or $\mathcal{M}^*$ the sets $\{M : M(f) \neq 0\}$ are open, because $M(f) - r = M(f - r)$; and because of the lattice homomorphism, $\{M : M(f) < r\}$ is the open set $\{M : M(f \wedge r) \neq r\}$.

On $\mathcal{M}^*$, we define the real-valued functions
$$\widehat{f}(M) = M(f) \qquad (f \in C^*,\ M \in \mathcal{M}^*).$$

The weak topology on $\mathcal{M}^*$ determined by these functions has for a base all sets of the form $\widehat{f}^{\leftarrow}(r, s)$; since $\widehat{f}^{\leftarrow}(r, s) = (\widehat{f} - t)^{\leftarrow}(-\epsilon, \epsilon)$, where $\epsilon = (s-r)/2$, $t = (s+r)/2$, the sets of the form $\widehat{g}^{\leftarrow}(-\epsilon, \epsilon)$ are themselves a base.

LEMMA.   *The weak topology on $\mathcal{M}^*$ coincides with the Stone topology and is compact.*

PROOF. $\widehat{f}^{\leftarrow}(-\epsilon, \epsilon)$ is the open set $\{M : |M(f)| < \epsilon\}$ in the Stone topology; and $\{M : M(f) \neq 0\}$ is the open set $\widehat{f}^{\leftarrow}(\mathbf{R} - \{0\})$ in the weak topology. So the topologies agree.

$\mathcal{M}^*$ is a Hausdorff space, because the continuous functions $\widehat{f}$ distinguish points. Since $C^*$ is a commutative ring with 1, every open cover of $\mathcal{M}^*$ has a finite subcover, as is well known. (A cover by sets $\{M : f_\alpha \notin M\}$ means that the $f_\alpha$ generate the whole ring (1), whence some finite subset generates (1).) Therefore $\mathcal{M}^*$ is compact.

**5.4   THEOREM**   *The mapping $f \mapsto \widehat{f}$ is an isomorphism of $C^*[\mathcal{S}]$ onto a dense subalgebra of $C^*(\mathcal{M}^*)$ which contains all constant functions.*[2]

PROOF. The mapping is obviously a homomorphism onto a subalgebra of $C^*(\mathcal{M}^*)$ containing the constants. If $f \neq 0$, then by Theorem 5.2, $\|\widehat{f}\| \neq 0$, whence $\widehat{f} \neq 0$; therefore the mapping is one-one. Since the subalgebra distinguishes points, it is dense, by the Stone-Weierstrass theorem.

**5.5   THEOREM**   *For each space $X$, there exists a compact, extremally disconnected space $K$—namely, $\mathfrak{M}(\overline{Q}^*(X))$—such that $\overline{Q}^*(X) \cong C(K)$.*

PROOF. Since $\overline{Q}^*(X)$ is metrically complete, it is isomorphic with $C^*(K)$, by the preceding theorem; and $K$ is compact (Lemma 5.3). Since $\overline{Q}^*(X)$ is Dedekind-complete (4.11), $K$ is extremally disconnected (4.12).

On taking rational completions, by the way, we get $\overline{Q}(X) \cong Q(K) = \overline{Q}(K)$ (Theorem 4.8) and hence $\overline{Q}^*(X) \cong C(K) = Q^*(K) = \overline{Q}^*(K)$. (Hence by 3.5, $K$ is extremally disconnected.)

---

[2]This statement clarifies the original.

## 5.6. The function $\check{f}$.

LEMMA. *Let $a \in C/M$, where $M \in \mathcal{M}$. If $|a|$ is not infinitely large, then there exists a unique real number $r$ such that $|a - r|$ is either infinitely small or zero.*

PROOF. By hypothesis, the set $\{s \in \mathbf{R} : s < a\}$ is nonvoid and has an upper bound in $\mathbf{R}$; let $r$ denote its supremum in $\mathbf{R}$. For each positive integer $n$,

$$r - \frac{1}{n} < a < r + \frac{1}{n},$$

whence $|a - r| < \frac{1}{n}$; therefore $|a - r|$ is infinitely small or 0. Uniqueness of $r$ is obvious.

COROLLARY. *For each $M \in \mathcal{M}$ and $f \in C^*$, there exists a unique real number $r$ such that $|M(f) - r|$ is either infinitely small or zero. (Cf. [Gillman & Jerison (1960), 7.6].)*

PROOF. There exists an integer $m$ such that $|f| \leq m$, whence $|M(f)| \leq m$. So $|M(f)|$ is not infinitely large.

For each $f \in C^*$, we now define a real-valued function $\check{f}$ on $\mathcal{M}$ as follows:

$$\check{f}(M) = r,$$

where $|M(f) - r|$ is infinitely small or zero.

REMARK. Observe that $\check{f}(M) < s$ implies $M(f) < s$ but that $M(f) < s$ implies only that $\check{f}(M) \leq s$.

## 5.7 Lemma

(1) *For each $M \in \mathcal{M}$,*
$$\tau^*(M) = \{f \in C^* : \check{f}(M) = 0\}$$
*is a maximal ideal in $C^*$.*

(2) *The weak topology on $\mathcal{M}$ (determined by the functions $\check{f}$) is Hausdorff.*

(3) *The mapping $\tau^* : \mathcal{M} \to \mathcal{M}^*$ defined in (1) is one-one and onto and satisfies $\check{f} = \widehat{f} \circ \tau^*$.*

PROOF. (1). If $|a|$ and $|b|$ are infinitely small (in $C/M$), so is $|a-b|$; and if $|c|$ is not infinitely large, then $|ac|$ is infinitely small. It follows that $\tau^*(M)$ is an ideal. Obviously, it is proper. To establish maximality, consider any $f \in C^*$ with $f \notin \tau^*(M)$. Then $f \notin M$, and so there exists $h \in C$ such that $M(hf) = 1$. By assumption, $|M(f)|$ is not infinitely small; hence $|M(h)|$ is not infinitely large, and so there exists a real number $r$ for which $|M(h) - r|$ is infinitely small or zero. Evidently, $r - 1 < M(h) < r + 1$. Hence for $h' = ((r-1) \vee h) \wedge (r+1)$, we have $h' \in C^*$ and $M(h') = M(h)$. Thus, $h'f - 1 \in M \cap C^* \subset \tau^*(M)$. This shows that $\tau^*(M)$ is maximal.

(2). If $M \neq M'$, there exist $f \in M \cap C^*$ and $g \in M' \cap C^*$ such that $f + g = 1$; then $\check{f}(M) = 0 \neq \check{f}(M')$. This shows that the functions $\check{f}$ distinguish points—whence the weak topology is Hausdorff—and that $\tau^*$ is one-one.

(3). We have just seen that $\tau^*$ is one-one. Given $M^* \in \mathcal{M}^*$, consider the set

$$L = \{(|f| - \epsilon) \vee 0 : f \in M^*, \epsilon > 0\}.$$

Now if $f_1, \cdots, f_n \in M^*$, then $f_1^2 + \cdots + f_n^2$ belongs to $M^*$, hence is not a unit of $C^*$, and so is not bounded away from 0. Given $\epsilon_1 > 0, \cdots, \epsilon_n > 0$, let $\epsilon = \min(\epsilon_1, \cdots, \epsilon_n)$; then

$$\bigcap_k Z((|f_k| - \epsilon_k) \vee 0) \supset \{x : f_1^2(x) + \cdots + f_n^2(x) < \epsilon^2\}$$

a nonvoid open set (in $\bigcap_k \mathrm{dom} f_k$). It follows that no linear combination in $C$ of the functions $(|f_k| - \epsilon_k) \vee 0$ can be equal to 1. The set $L$ therefore generates a proper ideal in $C$ and hence it is contained in some maximal ideal $M$.

For $f \in M^*$, we have $(|f| - \epsilon) \vee 0 \in M$, whence $|M(f)| \leq \epsilon$, for all $\epsilon > 0$. Hence $\check{f}(M) = 0$. Thus, $M^* \subset \tau^*(M)$l by maximality, $M^* = \tau^*(M)$. So $\tau^*$ is onto.

As we have just seen, $\widehat{f}(\tau^*(M)) = 0$ implies $\check{f}(M) = 0$. Hence, clearly, $\widehat{f}(\tau^*(M)) = r$ implies $\check{f}(M) = r$ for all $r$. Thus, $\check{f} = \widehat{f} \circ \tau^*$.

5.8 THEOREM *The weak topology on $\mathcal{M}$ coincides with the Stone topology and is compact; and the mapping $\tau^* : \mathcal{M} \to \mathcal{M}^*$ is a homeomorphism.*

PROOF. Since $\widehat{f}(-\epsilon, \epsilon) = \bigcup_{0 < \delta < \epsilon}\{M : |M(f)| < \delta\}$ (see 5.6, *Remark*) the Stone topology on $\mathcal{M}$ contains the Hausdorff weak topology. Since $C$ is a commutative ring with 1, $\mathcal{M}$ is compact in the Stone topology. Therefore the identity mapping from $\mathcal{M}$ in the Stone topology to $\mathcal{M}$ in the weak topology is a homeomorphism.

Because each function $\widehat{f} \cdot \tau^* (= \check{f})$ is continuous and the functions $\widehat{f}$ determine the topology of $\mathcal{M}^*$ (5.3), $\tau^*$ is continuous (see, e.g., [Gillman & Jerison (1960), 3.8]). Since $\tau^*$ is one-one from the compact space $\mathcal{M}$ to the Hausdorff space $\mathcal{M}^*$, it is a homeomorphism.

The proofs here of 5.7 and 5.8 give information about the behavior of functions. The homeomorphism between $\mathcal{M}$ and $\mathcal{M}^*$ will also follow from Theorem 10.18 below, whose proof has a more algebraic flavor. (To apply 10.18, note that $C^*$ is semi-simple (as is easily seen) and that $\mathcal{M}^*$ is Hausdorff (5.3).)

**5.9 Corollary** *The maximal ideal spaces of $Q(X)$ and $Q^*(X)$ are homeomorphic; so are those of $\overline{Q}(X)$ and $\overline{Q}^*(X)$.*

So are those of $C(X)$ and $C^*(X)$ as is well known [Gillman & Jerison (1960)].

We now consider a family $\mathcal{S}' \supset \mathcal{S}$ such that $C^*[\mathcal{S}]$ is metrically dense in $C^*[\mathcal{S}']$. With $C, C^*, \mathcal{M}$ and $\mathcal{M}^*$ as before, we also put $C' = C[\mathcal{S}'], \mathcal{M}' = \mathfrak{M}(C'), C'^* = C^*[\mathcal{S}']$, and $\mathcal{M}'^* = \mathfrak{M}(C'^*)$.

**5.10 Theorem** *If $C^*$ is metrically dense in $C'^*$, then $\mathcal{M}^*$ is homeomorphic with $\mathcal{M}'^*$, under the mapping $\sigma : M \mapsto \mathrm{cl}_{C'^*} M$. Hence $\mathcal{M}, \mathcal{M}', \mathcal{M}^*$, and $\mathcal{M}'^*$ are all homeomorphic.*

Proof. Because of continuity, $\mathrm{cl}\, M$ contains sums and products—i.e., it is an ideal. For maximality, we show that $C'^*/\mathrm{cl}\, M$ is the field $\mathbf{R}$. Given $f' \in C'^*$, choose $f_n \in C^*$ with $f_n \to f'$. By 5.1, $r_n = M(f_n) \in \mathbf{R}$. Since $|r_m - r_n| = M(|f_m - f_n|) \leq \|f_m - f_n\|$, $(r_n)$ is a Cauchy sequence in $\mathbf{R}$ and hence converges to some $r \in \mathbf{R}$. The sequence $(f_n - r_n)$ in $M$ then converges to $f' - r$ and so $f' - r \in \mathrm{cl} M$.

Thus, $\sigma$ maps $\mathcal{M}^*$ into $\mathcal{M}'^*$. Obviously it is one-one. We prove that $\sigma$ is onto. Given any $M' \in M'^*$, define $M = M' \cap C^*$. Evidently, $M$ is an ideal in $C^*$. For $f \in C^*$, we have $r = M'(f) \in \mathbf{R}$, by 5.1, whence $f - r \in M' \cap C^* = M$. Therefore, $C^*/M = \mathbf{R}$, and so $M$ is a maximal ideal in $C^*$. Now, we observe that $\mathrm{cl}\, M'$ is an ideal (again, by continuity) and proper (since the $\frac{1}{2}$-neighboorhood of 1 contains only units); by maximality $\mathrm{cl}\, M' = M'$. So we have $\mathrm{cl}\, M \subset \mathrm{cl}\, M' = M'$; since $\mathrm{cl}\, M$ is maximal, $\mathrm{cl} M = M'$. Thus, $\sigma$ is onto.

Given $f' \in C'^*$ and $\epsilon > 0$, let $f \in C^*$ satisfy $\|f' - f\| < \epsilon$; then on $\mathcal{M}'^*$, $|\widehat{f}(M')| < \epsilon$ implies $|\widehat{f}'(M')| < 2\epsilon$. Therefore the functions $\widehat{f}$, for $f \in C^*$ alone, determine the topology of $\mathcal{M}'^*$. Since $M(f) = (\mathrm{cl}\, M)(f)$, $\sigma$ carries the basic open set $\{M : |\widehat{f}(M)| < \epsilon\}$ in $\mathcal{M}^*$ onto the basic open set $\{\mathrm{cl}\, M : |\widehat{f}(\mathrm{cl}\, M)| < \epsilon\}$ in $\mathcal{M}'^*$. Consequently, $\sigma$ is a homeomorphism.

**5.11 Corollary** *For any space $X$, the maximal ideal spaces of $Q(X)$, $Q^*(X)$, $\overline{Q}(X)$, and $\overline{Q}^*(X)$ are all homeomorphic.*

# 6. Maximal Ideal Spaces of Subalgebras of $C^*(X)$

**6.1.** z-DETERMINING SUBALGEBRAS. Let $A$ be a *subalgebra* of $C^*(X)$: we shall use this term to include the condition $\mathbf{R} \subset A$.

An ideal $I$ in $A$ for which $Z(I)$ is nonempty is said to be *fixed*. For each $x \in X$, the set
$$M_x = \{f \in A : f(x) = 0\}$$
is a (fixed) maximal ideal in $A$ whose residue field is $\mathbf{R}$; for the mapping $f \to f(x)$ ($f \in A$) is an epimorphism of $A$ onto $\mathbf{R}$. When $A$ distinguishes the points of $X$, the mapping $x \to M_x$ from $X$ to the fixed maximal ideals is one-one.

For $S \subset X$, we write
$$M_S = \{M_x : x \in S\}.$$
In particular, $M_X$ denotes the set of all fixed maximal ideals in $A$.

We shall say that the subalgebra $A$ is *z-determining* if the zero-sets of $A$ form a base for the closed sets in $X$, i.e., if every closed set is an intersection of zero-sets of $A$. In particular, $C^*(X)$ itself is z-determining. Any z-determining subalgebra determines the topology of $X$ (i.e., as its weak topology). On the other hand, the polynomials in $C^*([0,1])$, for example, determine the topology of $[0,1]$ but are not z-determining. Also, a z-determining subalgebra need not be dense. An example is the algebra of all $f \in C^*(X)$ for which $f(p) = f(q)$, where $p$ and $q$ are fixed distinct points of $\beta X - X$. (This follows from the Stone-Weierstrass Theorem.)

**6.2 LEMMA** *If $X$ is compact and if the subalgebra $A$ distinguishes points and is a sublattice of $C^*(A)$, then $A$ is z-determining.*

PROOF. Let $F$ be a closed set in $X$ and let $p \in X - F$. For each $x \in F$ there exists $f_x \in A$ such that $f_x(p) > 0$ and $f_x(x) < 0$. Define $g_x = f_x \vee 0$. Then $g_x \in A$, $g_x(p) > 0$, and $Z(g_x)$ is a neighborhood of $x$. A finite number of these zero-sets—say $Z(g_{x_1}), \cdots, Z(g_{x_n})$—cover $F$. Then $Z(g_{x_1} \cdots g_{x_n})$ contains $F$ but not $p$. Thus, $F$ is an intersection of zero-sets of $A$. ∎

**6.3. THE STONE TOPOLOGY.** Recall that $\mathfrak{M}(A)$ denotes the set of all maximal ideals in $A$. For $f \in A$, we put
$$\Gamma(f) = \{M \in \mathfrak{M}(A) : f \notin M\}.$$
The sets $\Gamma(f)$ constitute a base for the Stone topology on $\mathfrak{M}(A)$.

LEMMA. *Let $A$ be a subalgebra of $C^*(X)$.*

1. *The subspace $M_X$ of fixed maximal ideals in $A$ is dense in $\mathfrak{M}(A)$.*

2. *If $X$ is compact, then $M_X = \bigcap_{Z(f)=\emptyset} \Gamma(f)$.*

PROOF. (1). If the basic open set $\Gamma(f)$ in $\mathfrak{M}(A)$ is not void, then $f \neq 0$; pick $x$ with $f(x) \neq 0$; then $M_x \in \Gamma(f)$.

(2). For each $x \in X$, $Z(f) = \emptyset$ implies $f \notin M_x$, i.e., $M_x \in \Gamma(f)$; so $M_x \in \bigcap_f \Gamma(f)$. Therefore $M_X \subset \bigcap_f \Gamma(f)$. Conversely, for each $M \in \bigcap_f \Gamma(f)$, $g \in M$ implies $Z(g) \neq \emptyset$; it follows easily that the family $\{Z(f) : f \in M\}$ has the finite intersection property. Since $X$ is compact, there exists $x \in Z(M)$. Then $M \subset M_x$ and therefore $M = M_x$. Consequently, $\bigcap_f \Gamma(f) \subset M_X$. ∎

**6.4. CONVEXITY.** We shall say that the subalgebra $A$ is *convex* if every element $\geq 1$ in $A$ is a unit of $A$. In case $X$ is pseudocompact, this is equivalent to: for $f \in A$, $Z(f) = \emptyset$ only if $f$ is a unit of $A$. Indeed, the latter assertion trivially implies convexity. Conversely, if $Z(f) = \emptyset$, then by pseudocompactness, $f^2 \geq r$ for some $r > 0$; then $f^2/r \geq 1$ and convexity implies that $f$ is a unit.

THEOREM. *If $X$ is compact, then the subalgebra $A$ of $C^*(X)$ is convex if and only if all the maximal ideals in $A$ are fixed.*

PROOF. If $A$ is convex, then $Z(f) = \emptyset$ implies that $f$ is a unit (since $X$ is compact), so that $\Gamma(f) = \mathfrak{M}(A)$; by the lemma, $M_X = \mathfrak{M}(A)$. Conversely, if $M_X = \mathfrak{M}(A)$, then by the lemma, $Z(f) = \emptyset$ implies that $f$ belongs to no maximal ideal, whence $f$ is a unit. ∎

COROLLARY. *If $X$ is compact, then $A$ is convex if and only if every proper ideal in $A$ is fixed.*

PROOF. Under convexity, a proper ideal is contained in a maximal ideal $M_x$. This yields the necessity; the converse is trivial. ∎

**6.5 Theorem**   *The subalgebra $A$ of $C^*(X)$ is $z$-determining if and only if $X$ is homeomorphic with $M_X$ under the mapping $x \to M_x$. Hence if $X$ is compact, $A$ is $z$-determining, and $\mathfrak{M}(A)$ is Hausdorff, then the homeomorphsim carries $X$ onto all of $\mathfrak{M}(A)$.*

Proof. If $A$ is $z$-determining, then the family of all sets $Z(f)$ for $f \in A$, is a base for the closed sets in $X$. Since $x \in Z(f)$ if and only if $f \in M_x$, the one-one mapping $x \to M_x$ carries this family onto the family of basic closed sets $\{M_x : f \in M_x\}$ ($f \in A$) in $M_X$. Therefore it is a homeomorphism.

Conversely, consider any closed set $F$ in $X$ and any point $p \in X - F$. If the mapping $x \to M_x$ is a homeomorphism, then $M_p \notin M_F$ and $M_F$ is closed in $M_X$. Since $M_p \in M_X$, this implies that $M_p \notin \mathrm{cl}_{\mathfrak{M}(A)} M_F$. Let $\Gamma(f)$ be a basic neighborhood (in $\mathfrak{M}(A)$) of $M_p$ disjoint from $M_F$; then $Z(f)$ contains $F$ but not $p$. Thus, $A$ is $z$-determining.

Under the additional conditions, the dense set $M_X$ in the Hausdorff space $\mathfrak{M}(A)$ is compact, hence closed, and hence all of $\mathfrak{M}(A)$. ∎

**6.6. Examples.**   Let $X = [0,1]$.

1. If $A$ is the algebra of all differentiable functions on $X$, then $A$ is convex and $z$-determining (but is not a lattice). Hence $\mathfrak{M}(A) \approx X$.

2. If $A$ is the algebra of all piecewise polynomials on $X$, then $A$ is $z$-determining but not convex. Hence not every maximal ideal is fixed. Indeed, the polynomial $x - 2$, for instance, has no zeros and so belongs to no fixed maximal ideal, but it is not a unit.

3. If $A$ is the algebra of all rational functions on $X$, then $A$ is convex, but not $z$-determining.

**6.7. Inverse limits.**   Let $\mathcal{S} = (S_a)$ be a filter base of dense sets in $X$. The index set is directed, $a \leq b$ meaning $S_a \supset S_b$. The direct limit

$$C^*[\mathcal{S}] = \varinjlim_a C^*(S_a)$$

was defined with respect to the monomorphism

(1) $$f_a \to f_b = f_a \circ \pi_a^b \qquad (a \leq b)$$

of $C^*(S_a)$ into $C^*(S_b)$, where $\pi_a^b$ denotes the injection of $S_b$ into $S_a$. Now, $f_a$ and $f_b$ have continuations to $\beta S_a$ and $\beta S_b$, and $\pi_a^b$ has a continuation from $\beta S_b$ onto $\beta S_a$ (onto because $S_b$ is dense in $\beta S_a$). Denote the extension by the original symbols; (1) holds in the extension, because it holds on the dense set $S_b$. For like reason, $\pi_a^b \circ \pi_b^c = \pi_a^c$ for $a \leq b \leq c$. With respect to the extended mappings, we have

$$C^*[\mathcal{S}] = \varinjlim_a C^*(\beta S_a).$$

In addition we are invited to consider the inverse limit

$$K = \varprojlim_a \beta S_a.$$

By definition, $K$ is the subspace of $\prod_a \beta S_a$ consisting of all points $x = (x_a)$ for which $\pi_a^b(x_b) = x_a$ $(a \leq b)$. As is well known, the projections $\pi_a$ from $K$ to $\beta S_a$ satisfy $\pi_a = \pi_a^b \circ \pi_b$ $(a \leq b)$. Also, $K$, as the inverse limit of compact spaces, is compact; and in this case, because all the mappings $\pi_a^b$ are onto, each $\pi_a$ is onto.

The members of $C^*[\mathcal{S}]$ are the elements of $f = (f_a)$ in $\prod_a C^*(\beta S_a)$ for which (1) holds. Hence if $f \in C^*[\mathcal{S}]$ and $x \in K$, then $f_a(x_a)$ is independent of $a$, and so

$$f'(x) = f_a(x_a)$$

defines $f'$ as a real-valued function on $K$.

These considerations lead to a new proof of Theorem 5.5.

6.8 THEOREM  *The mapping $f \to f'$ is an isomorphism of $C^*[\mathcal{S}]$ onto a dense subalgebra and sublattice of $C(K)$. Consequently, this subalgebra is $z$-determining, and $\mathfrak{M}(C^*[\mathcal{S}]) \approx K$.*

PROOF. Clearly, the mapping $f \to f'$ is a ring homomorphism. Fix $a$; since $\pi_a$ is onto, the mapping is one-one. Also, the equation $f' = f_a \circ \pi_a$ holds for all $f$; therefore $f'$ is the composition of two continuous functions and hence is continuous. Evidently, $A = \{f' : f \in C^*[\mathcal{S}]\}$ is a subalgebra of $C(K)$. If $x \neq y$ in $K$, then $x_a \neq y_a$ for some $a$ and then $f_a(x_a) \neq f_a(y_a)$ for some $f_a \in C^*(\beta S_a)$. Therefore $A$ distinguishes points; by the Stone-Weierstrass Theorem, $A$ is dense. Clearly, $f' \vee g' = (f_a \vee g_a) \circ \pi_a = (f \vee g)' \in A$; therefore $A$ is a sublattice of $C(K)$. By Lemma 6.2, $A$ is $z$-determining. Since $\mathfrak{M}(A)$ is a Hausdorff space (5.3), $\mathfrak{M}(A) \approx K$, by Theorem 6.5.

REMARK.  More generally, if $(K_a)$ is any inverse system of compact spaces with respect to onto maps $\pi_a^b$, then $(C(K_a))$ is a direct system of rings with respect to the monomorphisms (1), and $B = \varinjlim_a C(K_a)$ is isomorphic with a dense subalgebra and sublattice of $C(K)$, so that $B$ is $z$-determining and $\mathfrak{M}(B) \approx K$.

6.9 THEOREM  *For each space $X$, the compact space*

$$K = \varprojlim_{S \in \mathcal{G}_0(\beta X)} \beta S$$

*($\mathcal{G}_0(\beta X)$ denoting the family of all dense $\mathcal{G}_\delta$'s in $\beta X$) is homeomorphic with $\varprojlim_{V \in \mathcal{V}_0(\beta X)} \beta V$, with $\varprojlim_{V \in \mathcal{V}_0(X)} \beta V$, with $\mathfrak{M}(Q^*(X))$, and with $\mathfrak{M}(\overline{Q}^*(X))$; also $K$ is extremally disconnected. Moreover, $\overline{Q}^*(X) \cong C(K)$; thus*

$$\varinjlim_{S \in \mathcal{G}_0(\beta X)} C(\beta S) \cong C(\varprojlim_{S \in \mathcal{G}_0(\beta X)} \beta S).$$

*Finally $\overline{Q}(X) \cong Q(K)$.*

PROOF. Take $C^*[\mathcal{S}]$ in Theorem 6.8 as $Q^*(X)$ and as $\overline{Q}^*(X)$ (see 4.6). The final assertion and the fact that $K$ is extremally disconnected follow as in 5.5.

A direct proof that $\varprojlim_{V \in \mathcal{V}_0(\beta X)} \beta V$ is extremally disconnected appears in E. C. Weinberg's thesis [Weinberg (1961)]. ∎

# 7. Invariant Norms

**7.1.   Hemi-norms.** Let $A$ be an *algebra*, by which we mean a commutative algebra over **R** with unity element 1. By a *hemi-norm* on $A$, we shall mean a mapping $\nu\colon A \to [0,\infty]$ such that for all $a,b \in A$ and all $r \in \mathbf{R}$,

($\alpha$)  $\nu(ra) = |r| \cdot \nu(a)$ \hfill ($r \neq 0$ or $\nu(a) < \infty$),

($\beta$)  $\nu(a+b) \leq \nu(a) + \nu(b)$,

($\gamma$)  $\nu(ab) \leq \nu(a) \cdot \nu(b)$ \hfill (unless $\nu(a) \cdot \nu(b)$ is undefined)[3],

($\delta$)  $\nu(1) = 1$,

($\epsilon$)  $\nu(a) = 0$ only if $a = 0$.

By ($\alpha$) and ($\delta$), $\nu(0) = 0$. The restriction of $\nu$ to the subalgebra of elements of finite hemi-norm is evidently a *norm* on that subalgebra.

Remarks. (i). We could assume only $\nu(1) < \infty$, and then achieve ($\delta$) by defining a new hemi-norm $\nu'$ by: $\nu'(a) = \sup_{b \neq 0} \nu(ab)/\nu(b)$. (ii). More generally, $A$ could be any commutative ring with 1, condition ($\alpha$) then applying only for $r$ an integer.

**7.2.   Invariant norms.**   Let $\nu$ be a norm on $A$. For $D \subset A$, we define

$$\nu_D(a) = \sup_{0 \neq d \in D} \frac{\nu(ad)}{\nu(d)} \qquad (a \in A).$$

By ($\gamma$), $\nu_D(a) \leq \nu(a)$.

We say that the norm $\nu$ is *invariant* if $\nu_D = \nu$ whenever $D$ is a dense ideal in $A$.

More generally, consider an algebra $B \supset A$; For $D \subset A$ and $b \in B$ such that $bD \subset A$, we define

$$\nu_D(b) = \sup_{0 \neq d \in D} \frac{\nu(bd)}{\nu(d)}.$$

---

[3]This condition was expressed somewhat differently in the original.

**7.3 Lemma**   *Let $\nu$ be an invariant norm on $A$, let $B$ be a ring of quotients of $A$, and let $b \in B$. Then*
$$\nu_D(b) = \nu_{b^{-1}A}(b)$$
*for each dense ideal $D$ in $A$ contained in $b^{-1}A$.*

**Proof.** Obviously, $\nu_D(b) \leq \nu_{b^{-1}A}(b)$. In the other direction, consider any positive $r < 1$. There exists $0 \neq d \in b^{-1}A$ such that
$$r \cdot \nu_{b^{-1}A}(b) \leq \frac{\nu(bd)}{\nu(d)}.$$

Since $\nu$ is invariant, $\nu(bd) = \nu_D(bd)$, and so there exists $0 \neq d' \in D$ such that
$$r \cdot \nu(bd) \leq \frac{\nu(bdd')}{\nu(d')}.$$

Since $bd' \in bD \subset A$, $\nu(bd')$ is defined and satisfies
$$\nu(bd') \geq \frac{\nu(bdd')}{\nu(d)};$$

consequently,
$$\nu_D(b) \geq \frac{\nu(bd')}{\nu(d')} \geq \frac{\nu(bdd')}{\nu(d') \cdot \nu(d)} \geq r \cdot \frac{\nu(bd)}{\nu(d)} \geq r^2 \cdot \nu_{b^{-1}A}(b).$$

This implies that $\nu_D(b) \geq \nu_{b^{-1}A}(b)$. ∎

**7.4 Theorem**   *An invariant norm $\nu$ on $A$ has a canonical extension to a hemi-norm $\overline{\nu}$ on $\mathcal{Q}(A)$, defined by*
$$\overline{\nu}(b) = \nu_{b^{-1}A}(b) \qquad (b \in \mathcal{Q}(A)).$$

**Proof.** For $a \in A$, we have $1 \in A = a^{-1}A$; this implies that $\overline{\nu}(a) \geq \nu(a)$ and hence that $\overline{\nu}(a) = \nu(a)$. Therefore, $\overline{\nu}$ is an extension of $\nu$.

Evidently, $\overline{\nu}$ satisfies conditions ($\alpha$) and ($\delta$) in the definition of hemi-norm.

($\beta$). Given $b_1, b_2 \in \mathcal{Q}(A)$, define $D = b_1^{-1}A \cap b_2^{-1}A$. Then $D$ is a dense ideal in $A$ contained in $(b_1 + b_2)^{-1}A$. Applying the lemma and the fact that property ($\beta$) holds for the norm $\nu$, we get
$$\overline{\nu}(b_1 + b_2) = \nu_D(b_1 + b_2) = \sup_{0 \neq d \in D} \frac{\nu((b_1 + b_2)d)}{\nu(d)}$$
$$\leq \sup_{0 \neq d \in D} \frac{\nu(b_1 d)}{\nu(d)} + \sup_{0 \neq d \in D} \frac{\nu(b_2 d)}{\nu(d)}$$
$$= \nu_D(b_1) + \nu_D(b_2) = \overline{\nu}(b_1) + \overline{\nu}(b_2).$$

($\gamma$). Given $b_1, b_2 \in \mathcal{Q}(A)$, define $D = (b_1 b_2)^{-1} A \cap b_2^{-1} A$. Then $D$ is a dense ideal in $A$ and $b_2 D \subset b_1^{-1} A$. Since $\bar{\nu}(0) = 0$, we may assume in the proof that $b_2 \neq 0$; then $b_2 D \neq 0$. Again applying the lemma, and the fact that $\nu(0) = 0$, we obtain

$$\bar{\nu}(b_1 b_2) = \nu_D(b_1 b_2) = \sup_{0 \neq d \in D} \frac{\nu(b_1 b_2 d)}{\nu(d)}$$

$$= \sup_{d \in D,\, b_2 d \neq 0} \frac{\nu(b_1 b_2 d)}{\nu(b_2 d)} \cdot \frac{\nu(b_2 d)}{\nu(d)}$$

$$\leq \sup_{0 \neq a \in b_1^{-1} A} \frac{\nu(b_1 a)}{a} \cdot \sup_{0 \neq d \in D} \frac{\nu(b_2 d)}{\nu(d)}$$

$$= \nu_{b_1^{-1} A}(b_1) \cdot \nu_D(b_2) = \bar{\nu}(b_1) \cdot \bar{\nu}(b_2).$$

($\epsilon$). If $\bar{\nu}(b) = 0$, then $\nu(bd) = 0$ for all $d \in b^{-1} A$. Thus, $b \cdot (b^{-1} A) = 0$. Since $\mathcal{Q}(A)$ is a ring of quotients of $A$, $b = 0$. ∎

**7.5 THEOREM** *Let $\nu$ be an invariant norm on $A$ and define $\mathcal{Q}^*(A) = \{b \in \mathcal{Q}(A) : \bar{\nu}(b) < \infty\}$.*

(1) *If $A \subset B \subset \mathcal{Q}^*(A)$, then $\bar{\nu}|B$ is an invariant norm on $B$.*

(2) *$\mathcal{Q}^*(A)$ is the largest normed rational extension of $A$.*

PROOF.

(1) Let $E$ be any dense ideal in $B$ and consider any $b \in B$. By 1.4, $E \cap A$ is a dense ideal in $A$ and therefore $D = E \cap b^{-1} A$ is a dense ideal in $A$. Since $\bar{\nu}$ is an extension of $\nu$, we see that $\nu_D(b) \leq \bar{\nu}_E(b)$. Using the lemma, we get

$$\bar{\nu}(b) = \nu_D(b) \leq \bar{\nu}_E(b) \leq \bar{\nu}(b).$$

Therefore, $\bar{\nu}|B$ is invariant.

(2) If $B'$ is a rational extension of $A$, then $A \subset B' \subset \mathcal{Q}(A)$. Let $\nu'$ be a norm on $B'$ that extends $\nu$. Then for any $b \in B'$,

$$\bar{\nu}(b) = \nu_{b^{-1} A}(b) = \nu'_{b^{-1} A}(b) \leq \nu'(b) < \infty.$$

Therefore, $b \in \mathcal{Q}^*(A)$.

∎

7.6  LEMMA   *If $A$ is a subalgebra of $C^*(X)$, $D$ an ideal in $A$, and $p \in \operatorname{coz} D$, then there is a function $d \in D$ such that*
$$d(p) = \|d\| = 1.$$

PROOF. Choose $f \in D$ such that $s = f(p) \neq 0$. Since $D$ is an ideal in the algebra $A$, $g = f^2/s^2 \in D$; evidently, $g \geq 0$ and $g(p) = 1$. Let $n$ be an integral upper bound for $g$. The real polynomial
$$P(r) = \frac{1}{n^n} r \cdot (n+1-r)^n$$
satisfies $P(r) \geq 0$ for $0 \leq r \leq n$ (in fact, $0 \leq r \leq n+1$). Also $\max_r P(r) = P(1) = 1$. Since $D$ is an ideal in $A$, the function
$$d = P(g)$$
in $C^*(X)$ belongs to $D$. For every $x \in X$, since $0 \leq g(x) \leq n$, we have $0 \leq d(x) \leq 1$; also, $d(p) = P(1) = 1$. Consequently, $\|d\| = 1$. ∎

7.7  LEMMA   *If $A$ is a subalgebra of $C^*(Y)$ and $D$ is a dense ideal in $A$ such that $\operatorname{coz} D$ is dense in $Y$, then for each $g \in C^*(Y)$,*
$$\|g\| = \|g\|_D.$$

PROOF. Since $\operatorname{coz} D$ is dense, $\|g\| = \sup_{y \in \operatorname{coz} D} |g(y)|$. Consequently, given $\epsilon > 0$, there exists $p \in \operatorname{coz} D$ such that $|g(p)| > \|g\| - \epsilon$. By Lemma 7.6, there exists $d \in D$ such that $d(p) = \|d\| = 1$. Hence
$$\|g\|_D \geq \frac{\|gd\|}{\|d\|} = \|gd\| \geq |g(p)d(p)| = |g(p)| \geq \|g\| - \epsilon.$$
This implies that $\|g\|_D \geq \|g\|$ and therefore that $\|g\|_D = \|g\|$. ∎

7.8  THEOREM   *Let $A$ be a subalgebra of $C^*(X)$. Then for the following propositions:*

(1) *$A$ is z-determining,*

(2) *$A$ has a z-determining ring of quotients $B$ in $C^*(X)$,*

(3) *Each dense ideal in $A$ has a dense cozero-set,*

(4) *The sup norm on $A$ is invariant,*

(5) *$A$ is metrically dense in $C^*(X)$,*

*we have $(1) \Rightarrow (2) \Rightarrow (3) \Rightarrow (4)$, and $(4 \text{ and } 5) \Rightarrow (3)$.*

PROOF. (1) *implies* (2). Trivial.

(2) *implies* (3). Let $D$ be a dense ideal in $A$. Consider any nonvoid basic open set $\operatorname{coz} g$, where $g \in B$. Since $D$ is dense in $B$, we have $gD \neq 0$; therefore $\operatorname{coz} g$ meets $\operatorname{coz} D$. Thus, $\operatorname{coz} D$ is dense.

(3) *implies* (4). This is an immediate consequence of Lemma 7.7.

(4 *and* 5) *implies* (3). Let $D$ be an ideal in $A$ for which $\operatorname{coz} D$ is not dense in $X$; then there exists $g \in C^*(X)$ such that $\|g\| = 1$ and $g(\operatorname{coz} D) = 0$. Since $A$ is metrically dense, there exists $f \in A$ with $\|g - f\| < 1/2$; then $\|f\| > 1/2$, while $|f(x)| < 1/2$ for $x \in \operatorname{coz} D$. Hence for $0 \neq d \in D$, we have (noting that $d$ vanishes outside of $\operatorname{coz} D$),

$$\frac{\|fd\|}{\|d\|} = \frac{1}{\|d\|} \sup_{x \in X} |f(x)d(x)| = \frac{1}{\|d\|} \sup_{x \in \operatorname{coz} D} |f(x)d(x)|$$
$$\leq \frac{1}{\|d\|} \sup_{x \in \operatorname{coz} D} |f(x)| \cdot \|d(x)\| \leq 1/2.$$

Consequently, $\|f\|_D = \sup_{0 \neq d \in D} \frac{\|fd\|}{\|d\|} \leq 1/2 < \|f\|$; thus, the sup norm is not invariant. ∎

**7.9. COUNTEREXAMPLES.** Let $X = [0, 1]$.

(2) *does not imply* (1). Let $A$ be the algebra of all functions $f \in C^*(X)$ for which $f(0) = f(1)$. Obviously, $A$ is not z-determining. (It does not even distinguish points.) However, $C^*(X)$ itself is a ring of quotients of $A$. For, if $0 \neq g \in C^*(X)$, then there exists $p \in (0, 1)$ such that $g(p) \neq 0$. Pick $f \in C^*(X)$ such that $f(p) \neq 0$ while $f(0) = f(1) = 0$; then $f \in A$ and $0 \neq gf \in A$.

(3) *does not imply* (2). Take $A = \mathbf{R}$.

(5) *does not imply* (3); *hence* (5) *does not imply* (4). Let $d$ be a function in $C^*(X)$ such that $d^{\leftarrow}(r)$ contains an interval $I_r$ for infinitely many $r$ including $r = 0$; and let $A$ be the algebra generated by $d$ and the polynomials. Then $A$ is metrically dense in $C^*(X)$. Now, obviously, the principal ideal $(d)$ in $A$ does not have a dense cozero-set; but we shall show that $(d)$ is a dense ideal. Consider any $g \in A$ for which $gd = 0$. Now $g$ has the form

$$g = \sum_{k=0}^{n} p_k d^k,$$

where the $p_k$ are polynomials. Since $0 = gd = \sum_k p_k d^{k+1}$, we have $\sum_k p_k(x) r^{k+1} = 0$ for all $x \in I_r$ and hence for all $x$. For each $x$, then, the polynomial $q_x(y) = \sum_k p_k(x) y^{k+1}$ has infinitely many zeros; so all its coefficients $p_k(x)$ are zero. Thus each $p_k = 0$ and therefore $g = 0$. This shows that $(d)$ is a dense ideal.

**7.10 Corollary** *(i) The sup norm on $C^*(X)$ is invariant;*

*(ii) its canonical extension to $Q^*(X)$ coincides with the sup norm on $Q^*(X)$ and*

*(iii) is invariant.*

*(iv) $Q^*(X)$ is the largest normed rational extension of $C^*(X)$.*

*(v) The ring $\overline{Q}^*(X)$, with the sup norm, has no proper normed rational extension.*

PROOF. Since dense ideals in $C^*(X)$ have dense cozero-sets, (i) is given in Theorem 7.8 (with $A = C^*(X)$). Assertions (iii) and (iv) now follow from 7.5.

To prove (ii), consider any $g \in Q^*(X)$; then $g \in C^*(V)$ for some $V$ dense and open in $X$. Since $D = g^{-1}C^*(X)$ is a dense ideal in $C^*(X)$, Lemma 7.7 (with $A = C^*(X), Y = V$) yields $\|g\|_D = \|g\|$, q.e.d.

To establish (v), we recall that $\overline{Q}^*(X) \cong C(K)$ for suitable $K$ (5.5 or 6.9). Hence by (i), the sup norm on $\overline{Q}^*(X)$ is invariant. (Recall that by 5.2, the isomorphism preserves the sup norm.) By Theorem 4.8, $\mathcal{Q}^*(\overline{Q}^*(X)) = \overline{Q}^*(X)$; therefore $\overline{Q}^*(X)$ has no proper normed rational extension (7.5). ∎

**7.11 Theorem** *If $A$ is a z-determining subalgebra of $C^*(X)$, then $Q_L(X) \subset \mathcal{Q}(A) \subset Q(X)$. Consequently, $\mathcal{Q}(A)$ is dense in $Q(X)$ and $\overline{Q}(X) = \mathcal{Q}(A) + \overline{Q}^*(X)$.*

PROOF. Consider any $g \in Q_L(X)$ (see 4.3), and define $D = g^{-1}A$. Given $x \in \mathrm{dom}\, g$, put $r = g(x)$. Since $g^{\leftarrow}(r)$ is open and $A$ is z-determining, there exists $d \in A$ such that $x \in \mathrm{coz}\, d \subset g^{\leftarrow}(r)$. Then $gd$ extends to $r \cdot d \in A$. Therefore $d \in D$ and so $x \in \mathrm{coz}\, D$. This shows that $\mathrm{coz}\, D \supset \mathrm{dom}\, g$; therefore $\mathrm{coz}\, D$ is dense. Consequently, $D$ is a dense ideal in $A$ and we have $g \in \mathrm{Hom}\, D \subset \mathcal{Q}(A)$. Thus $Q_L(X) \subset \mathcal{Q}(A)$.

By Theorem 7.8, dense ideals in $A$ have dense cozero-sets. Therefore $\mathcal{Q}(A) \subset Q(X)$ (see 2.7). The remaining assertions follow from 4.7. ∎

REMARK. It is easily seen that the embedding of $\mathcal{Q}(A)$ in $Q(X)$ of 2.7 preserves the sup norm.

**7.12. Another example of a non-invariant norm.** Let $B$ be the Banach algebra of all continuously differentiable real functions on $[0, 1]$, with norm

$$\nu(f) = \|f\| + \|f'\|$$

(prime denoting derivative). Let $D$ be the principal ideal generated by the function $i(x) = x$; then every function in $D$ vanishes at 0. Obviously, $D$ is dense.

Given $0 \neq d \in D$, we have $r = \|d\| > 0$. If we look at a point $x$ for which $|d(x)| = r$, we see from the mean-value theorem that $\|d'\| \geq r/x \geq r$; consequently, $\nu(d) \geq 2r$. Next,

$$\nu(id) = \|id\| + \|(id)'\| \leq \|id\| + \|id'\| + \|d\| \leq r + \nu(d),$$

so that
$$\frac{\nu(id)}{\nu(d)} \leq \frac{r}{\nu(d)} + 1 \leq 3/2.$$
Therefore, $\nu_D(i) \leq 3/2 < 2 = \nu(i)$. Thus, $\nu$ is not invariant. ∎

# 8. Quasi-Real Rings

**8.1. FORMALLY REAL RINGS.** We shall say that a ring $A$ (commutative, with 1) is *formally real* if for all $a_1, \ldots, a_n \in A$, $\sum_k a_k^2 = 0$ implies all $a_k = 0$. Obviously, a formally real ring contains no nonzero nilpotent elements, i.e., is semi-prime.

THEOREM.   *If $A$ is formally real, then $\mathfrak{Q}(A)$ is formally real.*

PROOF. Let $\sum_k b_k^2 = 0$, where each $b_k \in \mathfrak{Q}(A)$. The ideal $D = \bigcap_k b_k^{-1} A$ in $A$ is dense. For $d \in D$,
$$\sum_k (b_k d)^2 = (\sum_k b_k^2) d^2 = 0 \; ;$$
hence each $b_k d = 0$. Thus, $b_k D = 0$, and therefore $b_k = 0$. ∎

**8.2. QUASI-REAL RINGS.**   Let $A$ be a partially ordered ring. We recall that the ordering is determined by its positive cone $P$, a subset of $A$ with the characteristic properties: $P + P \subset P$, $PP \subset P$, $P \cap -P = \{0\}$.

We shall say that the order on $A$ is *quasi-real* provided that all squares are positive, i.e., $a^2 \in P$ for all $a$; we also refer to $A$ itself, or to $P$, as quasi-real. Obviously, any *total* order is quasi-real: given $a$, we have $\pm a \geq 0$ and hence $a^2 \geq 0$.

If a ring $A$ admits quasi-real orderings, then it admits a smallest one: the intersection of them all. On the other hand, the union of a chain of quasi-real orders is, clearly, quasi-real; hence by Zorn's lemma, every quasi-real order is contained in a maximal one.

The set of all elements expressible as a sum of squares is denoted by $P_0$.

**8.3 THEOREM**   *A formally real ring $A$ admits a quasi-real ordering; for the smallest such, the positive cone is the set $P_0$. Conversely, any quasi-real semi-prime ring is formally real.*

PROOF. Clearly, $0 \in P_0$, $P_0 + P_0 \subset P_0$, and $P_0 P_0 \subset P_0$. Finally, if $a \in P_0 \cap -P_0$, then, since $A$ is formally real, $a = 0$. Thus, $P_0$ is a positive cone. Obviously, it is quasi-real and, in fact the smallest such.

Conversely, $\sum_k a_k^2 = 0$ implies $0 \leq a_i^2 \leq \sum_k a_k^2 = 0$, so that $a_i^2 = 0$, whence $a_i = 0$. ∎

We remark that a quasi-real ring need not be semi-prime; 9.2 below contains an example.

8.4  THEOREM   *Let $A$ be quasi-real and let $a \in A$. A necessary and sufficient condition that the order on $A$ be extendable to one in which $a$ is positive is that*

  (1)   *$P_a \cap -P_a = \{0\}$, where $P_a = \{pa + q : p \geq 0, q \geq 0\}$; and $P_a$ is the positive cone of the smallest such extension.*

  *A sufficient condition that the order be so extendable is*

  (2)   *for $x \in A$, $x^2 a \leq 0$ implies $x = 0$;*

*when $A$ is a regular ring, a sufficient condition is*

  (3)   *for $e$ an idempotent, $ea \leq 0$ implies $e = 0$;*

*when $A$ is a field, a sufficient conditon is*

  (4)   *$a \not\leq 0$.*

PROOF. The positive cone of such an extension clearly contains $P_a$; therefore (1) is necessary. Conversely, $P_a + P_a \subset P_a$; and because $a^2 \geq 0$, $P_a P_a \subset P_a$. Hence if (1) holds, then $P_a$ is a positive cone, and evidently the smallest one containing $P$ and $a$.

Now assume (2), and consider any $b \in P_a \cap -P_a$. There exist $p, q, p', q' \geq 0$ such that $b = pa + q = -(p'a + q')$; then $(p+p')a = -(q+q') \leq 0$. Therefore, $(p+p')^2 a \leq 0$; by (2), $p + p' = 0$ and hence $q + q' = 0$. It follows that $p = q = 0$, so that $b = 0$. Thus (1) holds.

If $A$ is regular, then (by definition) given $x \in A$, there exist $y \in A$ such that $x^2 y = x$; then $xy$ is idempotent. Hence if $x^2 a \leq 0$, then $xya = x^2 y^2 a \leq 0$, so that (3) implies $xy = 0$ and therefore $x = x(xy) = 0$; thus, (2) holds. Finally, if $A$ is a field, then the only idempotents are 0 and 1, whence (4) implies (3). ■

We shall have more to say about regular rings in Chapter 10.

A number of known results about fields may be deduced easily from the theorem.

8.5  COROLLARY   *Any maximal quasi-real order on a field is total. Hence every quasi-real order on a field can be extended to a total order.*

8.6  COROLLARY   *(Artin-Schreier). Every formally real field can be totally ordered.*

PROOF. 8.3 and 8.5. (Cf. [Bourbaki (1952)].)

8.7  COROLLARY   *(Artin). In a formally real field, any element that is positive in every total order is a sum of squares.*

PROOF. Consider any $a \notin P_0$; then $-a \notin -P_0$. By (4), the quasi-real order determined by $P_0$ (8.2) extends to a maximal one in which $a$ is negative—that is, to a total order (8.5) in which $a$ ($\neq 0$) is not positive.

We remark that this proof is more elementary than Artin's [Artin (1926)] since it is carried out entirely within the field and does not involve the notion of real closure.

One may ask whether Artin's theorem generalizes to regular rings: in a formally real regular ring, is an element that is positive in every maximal quasi-real order necessarily a sum of squares? The answer is "no"; a counterexample is described in 10.8 below.

**8.8. CONVEXITY.** Let $A$ be a partially ordered ring. As in 6.8, we shall say that $A$ is *convex* provided that every element $\geq 1$ is a unit. When $1 \geq 0$, as when $A$ is quasi-real, this implies that $A$ contains the rationals.

Next, we say that an ideal $I$ is a *convex ideal* if $0 \leq a \leq b$ and $b \in I$ implies $a \in I$. As is well known (and easily seen), convexity of $I$ is the necessary and sufficient condition that $A/I$ be partially ordered in a natural way: the positive elements are the classes $I(x)$ for $x \geq 0$ in $A$.

**8.9 THEOREM** *A quasi-real ring is convex if and only if all its maximal ideals are convex ideals.*

PROOF. Assume the ring convex, consider any maximal ideal $M$, and let $0 \leq a \leq b$, with $a \notin M$. There exists $x \in A$ such that $1 - xa \in M$. Then $m = 1 - x^2a^2 \in M$ and we have $1 = x^2a^2 + m \leq x^2b^2 + m \notin M$. Therefore $b \notin M$.

Conversely, assume maximal ideals convex. If $a \geq 1 \geq 0$, then $a$ belongs to no maximal ideal and hence is a unit. ∎

**8.10. THE SUBRING $A^*$ OF BOUNDED ELEMENTS.** Let $A$ be ordered. An element $a$ is *bounded* if there exists a natural number $m$ such that $-m \leq a \leq m$. The set of all bounded elements of $A$ is denoted by $A^*$.

THEOREM. *If $1 \geq 0$ in $A$, then $A^*$ is a subring of $A$.*

PROOF. $0 \in A^*$; so $A^*$ is not empty. Let $a, b \in A^*$; then $-m \leq a \leq m$ and $-n \leq b \leq n$ for some natural $m$ and $n$. Then $-(m+n) \leq a - b \leq m + n$. Next, $0 \leq (m-a)(n+b)$; since $m \geq 0$ and $n \geq 0$ in $A$,

$$ab \leq mn + mb - an \leq 3mn \ .$$

Similarly, $0 \leq (m-a)(n-b)$, whence $ab \geq -3mn$. Therefore $A^*$ is a ring.

**8.11. INFINITESIMALS.** Let $A$ be ordered. An element $a \in A$ will be called *infinitesimal* if $-1 \leq na \leq 1$ for all natural $n$. Trivially, infinitesimals are bounded.

LEMMA. *Let $A$ be ordered and let $a \in A$. If there exist an element $e$ and an integer $n$ such that $ena \geq e > 0$, then $a$ is not an infinitesimal.*

PROOF. If $a$ is infinitesimal, then $1 \geq 2na$, whence $e \geq 2ena > ena \geq e$.

8.12 THEOREM  *Let $A$ be quasi-real. If*
 *(1) for all $a > 0$, there exists an idempotent $e \neq 0$ and an integer $n$ such that $na \geq e$, then every infinitesimal in $A$ is nilpotent.*

PROOF. If $b$ is infinitesimal, then $a = b^2$ is infinitesimal. If $a \neq 0$, then $a > 0$ and $na \geq e$ as in (1) implies $e = e^2 > 0$ and then $ena \geq e > 0$, contradicting the lemma.

8.13 THEOREM  *Let $A$ be quasi-real and convex, and let $S$ be a subring containing $A^*$. Then:*

(1) *$S$ is quasi-real and convex.*

(2) *$x^2 \in S$ implies $x \in S$.*

(3) *$A$ is a ring of quotients of $S$.*

(4) *All residue fields of $S$ can be totally ordered; and those of $A^*$ are embeddable in $\mathbf{R}$.*

(5) *All infinitesimals of $A$ lie in the radical of $A^*$.*

PROOF. (1). Trivially, $S$ is quasi-real. By convexity of $A$, if $x \geq 1$, then $x^{-1}$ exists in $A$; we prove:

(6) if $x \geq 1$, then $0 \leq x^{-1} \leq 1$, whence $x^{-1} \in A^*$.
In fact, $x^{-1} = x(x^{-1})^2 \geq 0$, whence $x^{-1} \leq xx^{-1} = 1$. It follows from (6) that $S$ is convex.

(2) By (6), $y = (x^2 + 1)^{-1}$ is a positive element of $S$; so is $2^{-1}$. Hence from the evident relation
$$-(x^2 + 1) \leq 2x \leq x^2 + 1 ,$$
we get $-1 \leq 2xy \leq 1$, so that $2xy \in A^* \subset S$, whence $xy \in X$. Therefore $x = (xy)(x^2 + 1) \in S$.

(3) Given $b \in A$ and $0 \neq b' \in A$, we are to find $a \in S$ such that $ba \in S$ and $b'a \neq 0$ (see 1.4). Define $a = (b^2 + 1)^{-1}$; as above, $a \in S$ and $ba \in S$. Finally, $b'a \neq 0$, because $b' \neq 0$ and $a$ is a unit.

(4) Since each maximal ideal is convex, its residue field is partially ordered as described above (8.8). Evidently, the order is quasi-real; by 8.5, it extends to a total order. In the case of $A^*$, the field is obviously archimedean and hence, as is well known, embeddable in $\mathbf{R}$.

(5) Let $a$ be infinitesimal, and consider any maximal ideal $M$ in $A^*$. By (4), $A^*/M \subset \mathbf{R}$. Evidently, $|M(a)| \leq 1/n$ for all natural $n$ ($> 0$). Hence $a \in M$. ∎

REMARK. The embedding referred to in (4) is in general not unique, because the partial order may have several extensions to total orders. Compare Corollary 9.5.

# 9. $\pi$-Rings

**9.1.** $\pi$-VALUES. Let $A$ be an ordered ring, and let $a \in A$; by a $\pi$-*value* of $a$, we shall mean any element $p \geq 0$ for which $p^2 = a^2$.

THEOREM. *If $A$ is quasi-real and semi-prime, then each element of $A$ has at most one $\pi$-value.*

PROOF. Suppose that $p^2 = q^2$, where $p \geq 0$ and $q \geq 0$. Define
$$s = p(p-q)^2, t = q(p-q)^2 ;$$
then $s + t = 0$. Since $A$ is quasi-real, $s \geq 0$ and $t \geq 0$; therefore $s = t = 0$. Hence $(p-q)^3 = s - t = 0$. Since $A$ is semi-prime, $p - q = 0$. ∎

**9.2.** $\pi$-RINGS. An ordered ring in which each element $a$ has a *unique* $\pi$-value $\pi a$ will be called a $\pi$-*ring*. Clearly, $\pi$-rings are quasi-real.

Trivially, in any $\pi$-ring:

(1) $(\pi a)^2 = a^2$,

(2) $\pi(\pi a + \pi b) = \pi a + \pi b$,

(3) $\pi[(\pi a)(\pi b)] = (\pi a)(\pi b)$,

(4) $\pi(-a) = \pi a$.

Also, $\pi a = a$ if and only if $a \geq 0$; in particular, $\pi 0 = 0$ and $\pi 1 = 1$.

Property (3) can be considerably strengthened. For we have $[\pi(ab)]^2 = (ab)^2 = [(\pi a)(\pi b)]^2$; by uniqueness,

(3′) $\pi(ab) = (\pi a)(\pi b)$.

EXAMPLE. It will be seen (Theorem 9.7) that a semi-simple convex $\pi$-ring satisfies $\pi a \geq a$ for all $a$. Here we exhibit a convex $\pi$-ring $S$ in which the inequality fails. $S$ will be the set of all polynomials $a + bx$, where $a$ and $b$ are rational, under the usual operations except that $x^2 = 0$; and the positive elements will be 0 and all $a + bx$ for which $a > 0$. Obviously, $S$ is not semi-prime. Each element of $S$ has a unique $\pi$-value: $\pi(a + bx) = (\operatorname{sgn} a)(a + bx)$. Each element $a + bx$ for which $a \neq 0$—in particular, each element $\geq 1$—has an inverse: $(a + bx)^{-1} = (1/a) - (b/a^2)x$; therefore $S$ is convex. Finally, we have $\pi x = 0 \not\geq x$.

**9.3 THEOREM** *Let $A$ be a ring in which $2x^2 = 0$ implies $x = 0$, and suppose that with each element $a \in A$ there is associated an element $\pi a \in A$ such that (1), (2), (3), and (4) hold. Then $A$ is formally real, and $A$ can be ordered to become a $\pi$-ring, with $\pi a$ the $\pi$-value of $a$.*

PROOF. Let $P = \{a \in A : \pi a = a\}$; we verify that $P$ is a positive cone. By (1), $(\pi 0)^2 = 0$, whence, by hypothesis, $\pi 0 = 0$; so $0 \in P$. By (2), $P + P \subset P$; and by (3), $PP \subset P$. Finally, consider any $a \in P \cap -P$; then $\pi a = a$ and $\pi(-a) = -a$. By (4), $a = -a$, that is, $2a = 0$; hence by hypothesis, $a = 0$. Consequently, $P \cap -P = \{0\}$. Thus, $P$ is a positive cone.

Since $\pi 0 = 0$, (2) yields $\pi(\pi a) = \pi a$; therefore each element $\pi a$ belongs to $P$ and hence by (1) is a $\pi$-value of $a$. Furthermore, $P$ is quasi-real.

Obviously, $A$ is semi-prime. By 9.1, $\pi$-values are unique; therefore $A$ is a $\pi$-ring. By 8.3, $A$ is formally real. ∎

**9.4. π-IDEALS.** An ideal $I$ in a $\pi$-ring is a *$\pi$-ideal* provided that $a \equiv b \pmod{I}$ implies $\pi a \equiv \pi b \pmod{I}$.

THEOREM. *Let $A$ be a $\pi$-ring. A prime ideal $J$ in $A$ is convex if and only if it is a $\pi$-ideal not containing 2; and when $J$ is convex, $A/J$ is totally ordered.*

PROOF. If $J$ is convex, then $2 \notin J$, since $2 \geq 1 \geq 0$. Suppose $a - b \in J$. Then $(\pi a)^2 - (\pi b)^2 = a^2 - b^2 \in J$; hence $\pi a = \pi b \in J$ or $\pi a - \pi b \in J$. In the first case, $\pi a \in J$ and $\pi b \in J$, by convexity; hence in either case, $\pi a - \pi b \in J$. So $J$ is a $\pi$-ideal.

Conversely, suppose that $J$ is a $\pi$-ideal not containing 2, and let $0 \leq a \leq b \in J$. Since $a - (a - b)$ belongs to the $\pi$-ideal $J$;

$$2a - b = a - (b - a) = \pi a - \pi(a - b) \in J \, .$$

Therefore, $2a \in J$, which implies $a \in J$. So $J$ is convex.

Finally, if $J$ is convex, then $A/J$ is partially ordered. Since $(a - \pi a)(a + \pi a) = 0 \in J$, either $a - \pi a \in J$ or $a + \pi a \in J$. Hence $J(a) = J(\pi a) \geq 0$ or $J(a) = J(-\pi a) \leq 0$. ∎

REMARK. The prime ideal (2) in the $\pi$-ring of integers is a $\pi$-ideal but is not convex. However, in most of the rings considered here, 2 is a unit.

**9.5 COROLLARY** *If $A$ is a convex $\pi$-ring and $S$ is a subring containing $A^*$, then:*

(1) *$S$ is a convex $\pi$-ring.*

(2) *All residue fields of $S$ are totally ordered; and those of $A^*$ appear canonically as (totally ordered) subfields of $\mathbf{R}$.*

PROOF. (1). By (1) of Theorem 8.13, $S$ is convex. By (2) of the same theorem, $a \in S$ implies $\pi a \in S$.

(2) The first statement is immediate from Theorem 9.4; the second follows upon observing again that the residue fields of $A^*$ are archimedean. ∎

**9.6. LATTICE-ORDERED RINGS.** A lattice-ordered ring is, by definition, a partially ordered ring that is also a lattice, i.e., in which $a \vee b$ (and hence $a \wedge b = -(-a \vee -b)$) exist for all $a$ and $b$. In such a ring, one defines $|a| = a \vee -a$; it satisfies $|a| \geq 0$ (for a proof see [Gillman & Jerison (1960), 5A]). The easily established law $(x+y) \vee (x+z) = x + (y \vee z)$ yields $2(a \vee b) = a + b + |a-b|$; in particular, $2(a \vee 0) = a + |a|$ and (dually) $2(a \wedge 0) = a - |a|$.

The next theorem gives a condition that a $\pi$-ring be lattice-ordered, with $|a| = \pi a$. (Cf. the example in 9.2.) We remark that a lattice-ordered ring need not be quasi-real and hence the equation $|a|^2 = a^2$ can fail. An example is the direct sum of **R** with itself, with $P = \{(x,y): x \geq y \geq 0\}$; here $(0,1)^2 = (0,1) \not\geq (0,0)$ (and $|(0,1)| = (2,1)$). In fact, we have:

(1) $\qquad\qquad\qquad |a|^2 = a^2$ if and only if $(a \vee 0)(a \wedge 0) = 0$.

(In the preceding example, $(0,1) \vee (0,0) = (1,1)$ and $(0,1) \wedge (0,0) = (-1,0)$.) To prove (1), note first that if $x \leq 0$ and $nx = 0$ for some positive integer $n$, then $x = 0$; for, the least such $n$ satisfies $-x = (n-1)x$, whence $-x \leq 0$. Applying this result to the identity $a^2 - |a|^2 = 4(a \vee 0)(a \wedge 0)$, noting that $(a \vee 0)(a \wedge 0) \leq 0$, we obtain (1).

9.7 THEOREM *If $A$ is a semi-simple convex $\pi$-ring, then $A$ is a lattice-ordered ring with $|a| = \pi a$ (and hence with $(a \vee 0)(a \wedge 0) = 0$).*

PROOF. By 8.9 and 9.4, all residue fields $A/M$ are totally ordered. Hence $M(x) = \pm M(\pi x)$. If $M(x) \geq 0$, then $M(x) = M(\pi x)$; and if this holds for all $M$, then by semi-simplicity, $x = \pi x \geq 0$. Since for each $M$, $M(\pi a \pm a) = M(\pi a) \pm M(a) \geq 0$, we get $\pi a \geq a$ and $\pi a \geq -a$.

Next, suppose that $u \geq a$ and $u \geq -a$. Then $M(u) \geq M(a)$ and $M(u) \geq -M(a)$, so that $M(u) \geq M(\pi a)$; as before, $u \geq \pi a$. This shows that $\pi a = a \vee -a$—that is, $\pi a = |a|$. Since 2 is a unit of $A$, $a \vee b = 2^{-1}(a + b + |a-b|)$ exists. Thus, $A$ is a lattice. ∎

**9.8. MAXIMAL IDEAL SPACES.** Recall that the sets $\Gamma(a) = \{M: M(a) \neq 0\}$ form a base for the Stone topology on $\mathfrak{M}(A)$. When each field $A/M \subset$ **R**, we may define a function $\widehat{a}: \mathfrak{M} \to$ **R** by $\widehat{a}(M) = M(a)$. Clearly, the mapping $a \mapsto \widehat{a}$ is a homomorphism from $A$ into $\mathbf{R}^{\mathfrak{M}(A)}$. Its kernel is the radical of $A$: hence the condition for this mapping to be a *monomorphism* is that $A$ be semi-simple.

THEOREM. *If $A^*$ is a convex $\pi$-ring, then $\widehat{a}$ is defined for $a \in A^*$ and is continuous— that is, the mapping $a \mapsto \widehat{a}$ is a canonical homomorphism (with kernel $\mathrm{rad}\, A^*$) of $A^*$ into $C^*(\mathfrak{M}(A^*))$.*

PROOF. $A^*$ contains the rationals; by 9.5, each field $A^*/M$ is totally ordered and $\subset \mathbf{R}$, whence $\widehat{a}$ is defined; clearly, for each rational $r$, $M(r) = r$. Since $A^*/M$ is totally ordered, we have for any $x \in A$, $M(x) = \pm M(\pi x)$; hence

$$(1) \qquad M(x) > 0 \text{ if and only if } M(\pi x + x) \neq 0.$$

To prove that $\widehat{a}$ is continuous, it suffices to show that $\widehat{a}^{\leftarrow}(r, s)$ is open for all rational $r$ and $s$. We have

$$\widehat{a}^{\leftarrow}(r, s) = \{M \colon M(a - r) > 0\} \cap \{M \colon M(s - a) > 0\};$$

and by (1), this is open. ∎

9.9  THEOREM  *The canonical injection of $A$ into $C(\mathfrak{M}(A))$, whenever it exists, has an extension to a canonical embedding of $\mathcal{Q}(A)$ into $\mathcal{Q}(\mathfrak{M}(A))$.*

PROOF. Identify $A$ with its image in $C(\mathfrak{M}(A))$. According to Corollary 2.7, we must show that for each dense ideal $D$ in $A$, $\cozD$ is dense in $\mathfrak{M}(A)$. Consider any basic open set

$$\Gamma(a) = \{M \colon M \in \coz a\} \neq \emptyset \,;$$

then $a \neq 0$. Since $D$ is dense, there exists $d \in D$ such that $ad \neq 0$. Since $A$ is semi-simple, there exists $M \in \mathfrak{M}(A)$ such that $ad \notin M$. Then $M \in \coz a \cap \coz d$, as required. ∎

# 10. Regular $\pi$-rings

**10.1. INTRODUCTION.** Recall that $A$ (commutative, with 1) is regular provided that for each $a$, there exists $b$ such that $a^2b = a$. The element $ab$ is idempotent; hence so is $1 - ab$. Since $a(1 - ab) = 0$, every element is either a zero-divisor or a unit, and every proper prime ideal is maximal. Regular rings are semi-simple, hence semi-prime. A *rationally complete* ring is regular if and only if it is semi-simple. (See 1.11.)

LEMMA. *If the elements $a$ and $b$ of a ring $A$ satisfy $a^2b = a$, then $b^2a$ is the unique element $c$ satisfying simultaneously $a^2c = a$ and $c^2a = c$, so that the idempotents $f = ac$ and $e = 1 - f$ satisfy $fa = a$, $fc = c$, and $ea = ec = 0$. Finally, the element $u = e + c$ is a unit, and $a^2u = a$.*

In particular [Gillman & Henriksen (1956)], if $A$ is regular, then for each element $a \in A$, there exists a unit $u$ such that $a^2u = a$.

PROOF. Direct substitution shows that $c = b^2a$ satisfies $a^2c = a$ and $c^2a = c$. If $a^2x = a$ and $x^2a = x$, then $x = ax^2 = a^2cx^2 = (a^2x)cx = (a^2c)cx = a^2xc^2 = ac^2 = c$, yielding the uniqueness. Clearly, $a^2u = a^2c = a$. Since $(e+c)(e+a) = e + f = 1$, $u = e + c$ is a unit. ∎

10.2 THEOREM *A regular quasi-real ring is convex.*

PROOF. Let $a \geq 1$. If an element $x$ satisfies $ax = 0$, then $0 \leq x^2 \leq ax^2 = 0$; hence $x^2 = 0$, and so $x = 0$. Thus $a$ is not a zero-divisor; therefore it is a unit. ∎

10.3 COROLLARY *If $A$ is a regular $\pi$-ring, then: $A$ is convex; $A$ is a lattice, with $|a| = \pi a$ and $(a \vee 0)(a \wedge 0) = 0$; any subring $S$ containing $A^*$ is a convex $\pi$-ring (hence has convex radical); $A$ is a ring of quotients of $S$ and the residue fields of $A^*$ are canonically embeddable in $\mathbf{R}$.*

PROOF. The nonparenthetical statements follow from 9.7, 9.5, and 8.13. Since intersections of convex ideals are convex, a convex $\pi$-ring has a convex radical (8.9). ∎

The question now arises when $A^*$ will be semi-simple—for then it can be represented as a subring of $C(\mathfrak{M}(A^*))$ (9.8). Accordingly, we now examine the radical of $A^*$.

10.4 THEOREM *If $A$ is a regular $\pi$-ring, then the radical of $A^*$ consists precisely of the infinitesimal elements (8.11) of $A$.*

PROOF. Let $J$ denote the radical of $A^*$. According to (5) of Theorem 8.13, all infinitesimals lie in $J$.

Conversely, consider any $a \in J$. Put $b = \pi a$; then $b^2 = a^2$, so that $b \in J$. For each natural $n$, $0 \leq nb \wedge 1 \leq nb \in J$; since $J$ is convex (Corollary 10.3), $nb \wedge 1 \in J$. Therefore $(nb - 1) \wedge 0 = (nb \wedge 1) - 1$ is a unit. But (see Corollary 10.3)

$$[(nb - 1) \wedge 0] \cdot [(nb - 1) \vee 0] = 0.$$

Hence $(nb - 1) \vee 0 = 0$—that is, $nb \leq 1$. This shows that $b$ is infinitesimal. Since $-b \leq a \leq b$, so is $a$. ∎

**10.5  COROLLARY**  *Let $A$ be a regular $\pi$-ring. If $A$ has no nonzero infinitesimals (for example, if (1) of 8.12 holds), then $A^*$ is canonically embedded in $C^*(\mathfrak{M}(A^*))$ and $A$ in $Q(\mathfrak{M}(A^*))$.*

PROOF. By the theorem, $A^*$ is semi-simple. The first result then follows from 9.8. Since $A$ is a rational extension of $A^*$ (10.3), we have $A \subset \mathcal{Q}(A^*)$, and the second result follows from 9.9. ∎

**10.6.  THE SET $E$ OF ALL IDEMPOTENTS.**  The set of all idempotents in a ring $A$ (not necessarily regular) will be denoted by $E$ or $E(A)$. This set may be ordered (whether or not $A$ is ordered), as follows: $e \leq f$ if and only if $ef = e$. Then $0 \leq e \leq 1$ and $E$ becomes a complemented distributive lattice, with $e \wedge f = ef$, $e \vee f = e + f - ef$, and $1 - e$ the complement of $e$. The corresponding Boolean ring has the operations $\dotplus$ and $\cdot$ defined by $e \dotplus f = e + f - 2ef$ and $e \cdot f = ef$.

When $A$ is a *quasi-real* ordered ring, the order on $E$ in $A$ coincides with the lattice order just described. For in $A$, we have $e = e^2 \geq 0$ and $1 - e = (1 - e)^2 \geq 0$, so that $0 \leq e \leq 1$. Hence $e \leq f$ implies $e = e^2 \leq ef \leq e \cdot 1 = e$; and, conversely, $ef = e$ implies $e \leq e + (f - e)^2 = f$.

**10.7  THEOREM**  *If $A$ is a regular ring with a maximal quasi-real order, and if the lattice $E$ of idempotents is complete, then $A$ is a $\pi$-ring.*

PROOF. Consider any $a \in A$. By Lemma 10.1, there is a unit $u$ such that $a^2 u = a$. Define

$$F = \{f \in E : fu \leq 0\}$$

and put

$$f_0 = \sup F \quad \text{in} \quad E.$$

Then $f \in F$ implies $f \leq f_0$ and hence $f f_0 = f$.

Next, define

(1) $$p = (1 - 2f_0)u$$

and $\pi a = a^2 p$. Then $p^2 = u^2$, so that $(\pi a)^2 = a^4 p^2 = a^4 u^2 = a^2$. It remains to show that $p \geq 0$, yielding $\pi a \geq 0$, so that $\pi a$ will be a $\pi$-value of $a$. Uniqueness of $\pi$-values is guaranteed by Theorem 9.1.

Consider any $e \in E$ with

(2) $$ep \leq 0;$$

we shall show that $e = 0$. Since the quasi-real order on $A$ is maximal, it will follow from (3) of Theorem 8.4 that $p \geq 0$, q.e.d.

For all $f \in F$, we have $efu \leq 0$. On the other hand, (1) yields $fp = -fu$, so that $efu = -efp \geq 0$, by (2). Therefore $efu = 0$. Since $u$ is a unit, $ef = 0$. Hence $f = (1-e)f \leq (1-e)f_0 \in E$. As this holds for all $f \in F$, we conclude that $f_0 \leq (1-e)f_0$; therefore $ef_0 = 0$. Now (1) and (2) yield $eu = ep \leq 0$. Thus, $e \in F$. So $e = ef_0 = 0$. ∎

**10.8. A COUNTEREXAMPLE.** The question was raised in 8.7 whether in a formally real regular ring, any element that is positive in every maximal quasi-real order is necessarily a sum of squares. We proposed the following counterexample. Let $F_n$ denote the field of rational functions, with real coefficients, in the $n$ indeterminates $x_1, \ldots, x_n$, and let $S$ be the direct sum of $F_1, F_2, \ldots$. It is easy to see that $S$ is formally real and regular, and that $E(S)$ is complete. By Theorem 10.7, $S$ is a $\pi$-ring in any maximal quasi-real order. Consider the element
$$a = (x_1^2, x_1^2 + x_2^2, \ldots).$$
Write $\pi a = (p_1, p_2, \ldots)$; since $(\pi a)^2 = a^2$, we have $p_n = \pm(x_1^2 + \cdots + x_n^2)$. Consider the idempotent $e_n$ whose $n^{\text{th}}$ component is 1 and all others 0; since $p_n = e_n(\pi a) \geq 0$, we must have $p_n = +(x_1^2 + \cdots + x_n^2)$. It follows that $\pi a = a$, so that $a > 0$.

Suppose now that $a$ were a sum of $n-1$ squares:
$$a = a_1^2 + \cdots a_{n-1}^2 \qquad (a_k \in S).$$
Let $f_k$ denote the $n^{\text{th}}$ component of $a_k$; then the functions $f_1, \ldots, f_{n-1}$ would have to satisfy
$$f_1^2 + \cdots f_{n-1}^2 = x_1^2 + \cdots + x_n^2.$$
It has since been shown that such a relation is not possible; Davenport [Davenport (1963)] established this in a special case, and later Cassels [Cassels (1964)] solved the general case.

**10.9 THEOREM** *If $A$ is rationally complete and semi-prime (hence regular), then $E(A)$ is complete.*

This will follow from Theorem 11.9. Assuming the result, we have:

**COROLLARY.** *A rationally complete semi-prime ring with maximal quasi-real order is a $\pi$-ring.*

We turn now to maximal ideal spaces.

**10.10  Theorem**   *Let A be regular and let S be a subring containing $E(A)$ and such that each proper prime ideal in S is contained in a unique maximal ideal in S. Then the maximal ideal spaces of A, S, and E are homeomorphic.*

Remark. If $e \in E$, then $e(1-e) = 0$. Hence any given prime ideal in $A$, $S$, or $E$ contains either $e$ or $1-e$. It follows that $\mathfrak{M}(E)$ and $\mathfrak{M}(A)$ are Hausdorff and totally disconnected. Since each ring contains 1, all three spaces have the Heine-Borel property (see end of 5.3). The result $\mathfrak{M}(A) \approx \mathfrak{M}(E)$ appears in [Morrison (1955)].

Proof. We shall prove that $\mathfrak{M}(S) \approx \mathfrak{M}(E)$. Since proper prime ideals in $A$ are maximal, the result for $S$ includes that for $A$.

Given any maximal ideal $M$ in $E$, let $MA$ denote the ideal generated by $M$ in $A$. This is a proper ideal. For if $1 \in MA$, then $1 = \sum_k e_k a_k$ for suitable $e_k \in M$ and $a_k \in A$; then $e = \prod_k (1-e_k) \notin M$, whereas $e = e \cdot 1 = 0$. Next $MA$ is a maximal (hence prime) ideal; for if $a = a^2 x \notin MA$, then $ax \notin MA$, whence $1 - ax \in MA$. It follows that $MA \cap S$ is a proper prime ideal in $S$ and hence by hypothesis is contained in a unique maximal ideal $M'$ in $S$. We have thus defined a mapping $M \to M'$—via

(1) $$M \to MA \cap S \to M'$$

—from $\mathfrak{M}(E)$ into $\mathfrak{M}(S)$.

This mapping is one-one. For, given two distinct maximal ideals in $E$, there exists $e$ such that one contains $e$ and the other $1 - e$, and no $M'$ contains both.

Next, we prove that the mapping is onto. Let $N$ be any maximal ideal in $S$. For $e \in E$, either $e$ or $1 - e$ belongs to $N$ and hence to $M = N \cap E$; therefore $M$ is maximal in $E$. Now consider any $s \in S$, with $s \notin N$. There exists $y \in A$ such that $s^2 y = s$. Then $sy \notin M$; therefore $1 - sy \in M \subset MA$, whence $s \notin MA$. We have shown that $MA \cap S \subset N$. By uniqueness, $N = M'$. This yields the result—and a byproduct as well:

(2) $$M' \cap E = M.$$

It remains to prove that the one-one mapping $M \to M'$ of $\mathfrak{M}(E)$ onto $\mathfrak{M}(S)$ is a homeomorphism. Since it carries the basic open set $\{M : e \notin M\}$ in $\mathfrak{M}(E)$ to the set

$$\Gamma(e) = \{M' : e \notin M'\} \qquad (e \in E)$$

in $\mathfrak{M}(S)$, we will complete the proof if we show that the sets $\Gamma(e)$ form a base in $\mathfrak{M}(S)$. Let $M' \in \mathfrak{M}(S)$ and let $\mathcal{B}$ be any open neighborhood of $M'$. For each $M'_\alpha \in \mathfrak{M}(S) - \mathcal{B}$, there exists $e_\alpha \in E$ such that $e_\alpha \in M'$ and $e_\alpha \notin M'_\alpha$. The open cover $(\Gamma(e_\alpha))$ of the quasi-compact set $\mathfrak{M}(S) - \mathcal{B}$ has a finite subcover, say $(\Gamma(e_k))$. Define $e = \prod_k (1 - e_k)$; then $M' \in \Gamma(e) \subset \mathcal{B}$. ∎

**10.11 LEMMA** *If $S$ is semi-simple and $\mathfrak{M}(S)$ is Hausdorff, then each proper prime ideal in $S$ is contained in a unique maximal ideal.*

PROOF. Distinct maximal ideals $M_1$ and $M_2$ have disjoint neighborhoods, say $\{M\colon s_1 \notin M\}$ and $\{M\colon s_2 \notin M\}$. Then $s_1 s_2 = 0$, so that no ideal contained in $M_1 \cap M_2$ is prime. (This result appears in [Gillman (1957)].) ∎

**10.12 THEOREM** *If $A$ is a regular $\pi$-ring and $A^*$ is semi-simple, then the maximal ideal spaces of $A$, $A^*$, and $E(A)$ are homeomorphic.*

PROOF. As we know, $E \subset A^*$ ($0 \le e \le 1$). By Corollary 10.3, $A^*$ is a convex $\pi$-ring. By Theorem 9.8, there is a canonical monomorphism $a \to \widehat{a}$ of $A^*$ into $C^*(\mathfrak{M}(A^*))$. Since the continuous functions $\widehat{a}$ distinguish points, $\mathfrak{M}(A^*)$ is a Hausdorff space. By Lemma 10.11, each proper prime ideal in $A^*$ is contained in a unique maximal ideal. Thus $S = A^*$ satisfies all the assumptions of Theorem 10.10. ∎

**10.13. RINGS OF CONTINUOUS FUNCTIONS.** We now apply some of the preceding results to $\mathrm{Q}(X)$ and its related rings. These rings come equipped with a natural order—the usual order, defined pointwise. In particular, $\mathrm{Q}$ and $\overline{\mathrm{Q}}$ are lattices and $\pi$-rings. In addition, the semi-prime rings $\mathrm{Q}$ and $\overline{\mathrm{Q}}$ are rationally complete (4.8) and hence regular (cf. 2.6, REMARK); by Theorem 10.9, the lattices $E(\mathrm{Q})$ and $E(\overline{\mathrm{Q}})$ (which, by what follows, coincide) are complete.

THEOREM. $E(\overline{\mathrm{Q}}(X)) = E(\mathrm{Q}(X))$.

PROOF. We must show that $E(\overline{\mathrm{Q}}) \subset E(\mathrm{Q})$. An idempotent of $\overline{\mathrm{Q}}$ is a continuous function $e$ from a dense $G_\delta$-set $S$ in $X$ into $\{0,1\}$. With each $s \in S$, associate an open neighborhood $V_s$ of $s$ in $X$ in which $e$ is constant. Consider the open sets

$$U_0 = \bigcup\{V_s \colon s \in S, e(s) = 0\}, \quad U_1 = \bigcup\{V_s \colon s \in S, e(s) = 1\}.$$

Since $U_0 \cup U_1$ contains $S$, it is dense. Now, evidently, $U_0 \cap U_1 \cap S = \varnothing$; since $S$ is dense, $U_0 \cap U_1 = \varnothing$. Therefore $e$ can be extended to an idempotent continuous on $U_0 \cup U_1$, hence belonging to $\mathrm{Q}$. ∎

**10.14 COROLLARY** $\mathfrak{M}(\overline{\mathrm{Q}}) \approx \mathfrak{M}(\mathrm{Q}) \approx \mathfrak{M}(E(\mathrm{Q}))$.

PROOF. 10.10 and 10.13.

As we see from (1) and (2) of 10.10, the homeomorphism from $\mathfrak{M}(\mathrm{Q})$ to $\mathfrak{M}(\overline{\mathrm{Q}})$ is given by $M \to \overline{M} = M\overline{\mathrm{Q}}$, its inverse by $\overline{M} \to M = \overline{M} \cap \mathrm{Q}$. ∎

10.15   THEOREM   *If a subring $A$ of $\overline{Q}(X)$ contains $E(Q)$ then the natural order $\geq$ on $A$ is maximal. In particular, the orders on $Q$ and $\overline{Q}$ are maximal.*

PROOF. Let $P$ denote the positive cone of some larger order, and consider any $a \in P$. Suppose that $a \not\geq 0$. Pick $x_0$ for which $a(x_0) < 0$. Let $U$ be a neighborhood of $x_0$ in $X$ such that $a(x) < 0$ for all $x \in U \cap \text{dom}\, a$. Define $e(x) = 1$ for $x \in U$, $e(x) = 0$ for $x \in X - \text{cl}\, U$. Then $e \in E(Q) \subset A$ and $ea \leq 0$; hence $ea \in -P$. But $e \geq 0$, so that $ea \in P$. Consequently $ea = 0$, and so $a(x_0) = 0$, a contradiction. ∎

10.16   COROLLARY   *Any regular ring $A$ between $E(Q)$ and $\overline{Q}$ is a sublattice of $\overline{Q}$; and if $A^*$ is semi-simple, then its maximal ideal space is homeomorphic to those of $A$ and $E(A)$ $(= E(Q))$.*

PROOF. According to Theorem 10.9, $E(A)$ is complete. Since the order on $A$ is maximal, $A$ is a π-ring, by Theorem 10.7. Hence $A$ is a sublattice of $\overline{Q}$. The rest is stated in Theorem 10.12. ∎

10.17   COROLLARY   *The maximal ideal spaces of $Q$, $Q^*$, $\overline{Q}$, $\overline{Q}^*$ and $E(Q)$ are all homeomorphic. (Cf. 5.11.)*

PROOF. We need only check that $Q^*$ and $\overline{Q}^*$ are semi-simple. Given $f \neq 0$, it is easy to find $g$ such that $g$ vanishes on a nonvoid open set and $f^2 + g^2$ is bounded away from 0—so that $g$ belongs to some maximal ideal $M$ while $f \notin M$. ∎

10.18   THEOREM   *Let $A$ be a convex π-ring satisfying $\pi a \geq a$ for all $a \in A$ (see 9.2, EXAMPLE). Let $S$ be a semi-simple subring containing $A^*$ and such that $\mathfrak{M}(S)$ is Hausdorff. Then $\mathfrak{M}(A)$ is homeomorphic with $\mathfrak{M}(S)$ under the mapping $\sigma$ defined by:*

$$\sigma(M) = \text{unique maximal ideal in } S \text{ containing } M \cap S \qquad (M \in \mathfrak{M}(A)).$$

PROOF. Since $M \cap S$ is a prime ideal in $S$, the mapping $\sigma$ is well defined (Lemma 10.11). To see that $\sigma$ is one-one, let $M \neq M'$ in $\mathfrak{M}(A)$. Choose $a, a' \in A$ such that $a \in M$, $a' \in M'$, and $a + a' = 1$. Then $\pi a \in M$, $\pi a' \in M'$, and $\pi a + \pi a' \geq 1$. Put $b = (\pi a + \pi a')^{-1}$ By (6) of Theorem 8.13, $b \geq 0$. Hence $0 \leq (\pi a) \cdot b \leq (\pi a + \pi a') \cdot b = 1$; thus $(\pi a) \cdot b \in A^* \subset S$, so that $(\pi a) \cdot b \in M \cap S \subset \sigma(M)$. Likewise, $(\pi a') \cdot b \in \sigma(M')$. Since $(\pi a) \cdot b + (\pi a') \cdot b = 1$, we have $\sigma(M) \neq \sigma(M')$.

Next, $\sigma$ is onto. For let $N \in \mathfrak{M}(S)$ be given. Define

$$N : A = \{a \in A : aA \subset N\}.$$

This is easily seen to be a prime ideal in $S$; and $N$ is the unique maximal ideal (in $S$) containing it. But $N : A$ is also a proper ideal in $A$ and hence is contained in some maximal ideal $M$ in $A$. Thus, $N : A \subset M \cap S \subset \sigma(M)$; by uniqueness, $\sigma(M) = N$. This shows that $\sigma$ is onto.

Next, we prove that $\sigma^{-1}$ is continuous. The problem reduces to this: given $a \in A$ and $M_0 \in \mathfrak{M}(A)$, with $a \notin M_0$, to find $s \in S$ for which $s \notin \sigma(M_0)$ and such that for all $M \in \mathfrak{M}(A)$, $a \in M$ implies $s \in \sigma(M)$. Choose $h \in A$ and $m \in M_0$ satisfying $m + ha = 1$; then $\pi m \in M_0$ and $\pi m + \pi(ha) \geq 1$. Put $b = (\pi m + \pi(ha))^{-1}$. Again $b \geq 0$, so that $0 \leq (\pi m) \cdot b \leq 1$, whence $(\pi m) \cdot b \in S$; also, $s \in S$, where

$$s = \pi(ha) \cdot b.$$

Since $(\pi m) \cdot b \in M_0 \cap S \subset \sigma(M_0)$ and $(\pi m) \cdot b + s = 1$, we have $s \notin \sigma(M_0)$. Next, by (3') of 9.2, $\pi(ha) = (\pi h)(\pi a)$. Therefore if $a \in M$, then $s \in M \cap S \subset \sigma(M)$.

Next, let $M \neq M'$ in $\mathfrak{M}(A)$. Then $\sigma(M) \neq \sigma(M')$. Since $\mathfrak{M}(S)$ is Hausdorff and $S$ is semi-simple, there exist $s, s' \in S$ such that $s \notin \sigma(M)$, $s' \notin \sigma(M')$, and $ss' = 0$. Then $s \notin M$, $s' \notin M'$, and $ss' = 0$. Therefore $\mathfrak{M}(A)$ is Hausdorff.

Finally, since $S$ is a commutative ring with 1, $\mathfrak{M}(S)$ is compact. Therefore the one-one, continuous mapping $\sigma^{-1}$ from $\mathfrak{M}(S)$ onto the Hausdorff space $\mathfrak{M}(A)$ is a homeomorphism. ∎

# 11. Boolean algebras

**11.1 Theorem** *If $A$ is Boolean then $\mathcal{Q}(A)$ is Boolean.*

Proof. Let $\phi \in \mathcal{Q}(A)$; then $\phi \in \operatorname{Hom} D$ for some dense ideal $D$ in $A$. For $d \in D$,

$$\phi^2(d) = \phi(\phi(d)) = \phi(\phi(d^2)) = \phi(d.\phi(d)) = \phi(d)^2 = \phi(d).$$

(This result appears in [Brainerd & Lambek (1959)] ). ∎

**11.2. Algebra of sets.** For any space $X$, we write

$$\mathcal{B}_0(X) = \text{Boolean algebra of all open-and-closed sets,}$$

with set-theoretic intersection and complement as Boolean meet and complement.

Next, recall that an open set is said to be *regular* if it is the interior of its closure. For $U \subset X$, let

$$U\check{\phantom{x}} = X - \operatorname{cl} U.$$

It is easy to see that a set $V$ is a regular open set if and only if $V = U\check{\phantom{x}}$ for some open $U$. Also, $U\check{\phantom{x}}\check{\phantom{x}} \supset U$, and $T \supset U$ implies $T\check{\phantom{x}} \subset U\check{\phantom{x}}$. It follows that $V$ is a regular open set if and only if $V\check{\phantom{x}}\check{\phantom{x}} = V$.

We put

$$\mathcal{B}(X) = \text{Boolean algebra of all regular open sets,}$$

with intersection as meet and $V \to V\check{\phantom{x}}$ as complementation. It is easily seen that $\mathcal{B}(X)$ is complete. In fact,

(1) $\qquad\qquad\qquad \mathcal{B}(X)$ is the completion of $\mathcal{B}_0(X)$,

as is well known [Halmos (1963), Theorem 11].

**11.3 Theorem** *If $A$ is semi-simple, then $E(A) \cong \mathcal{B}_0(\mathfrak{M}(A))$, under the mapping $\Gamma$ defined by*

$$\Gamma e = \{M \in \mathfrak{M}(A) : e \notin M\}.$$

*Hence the completion of $E(A)$ is $\mathcal{B}(\mathfrak{M}(A))$.*

PROOF. This theorem is known in greater generality [Jacobson (1956), p 209] but we outline a proof for convenience. Clearly, $\Gamma$ is a homomorphism into $\mathcal{B}_0$. By semi-simplicity, $\Gamma$ is one-one. Given $\mathfrak{B} \in \mathcal{B}_0$, put $I = \bigcap \mathfrak{B}$ and $J = \bigcap (\mathfrak{M} - \mathfrak{B})$. For $M \in \mathfrak{M}$, $M \supset I$ if and only if $M \in \mathfrak{B}$, and $M \supset J$ if and only if $M \notin \mathfrak{B}$. Hence $A$ is the direct sum of $I$ and $J$, whence $\mathfrak{B} = \Gamma e$ for some $e \in E$. So $\Gamma$ is onto. Thus, $\Gamma$ is an isomorphism. The final conclusion follows from (1). ∎

**11.4 COROLLARY** *If $A$ is a semi-simple ring and if $\mathfrak{M}(A)$ is a totally disconnected Hausdorff space, then $\mathfrak{M}(A) \approx \mathfrak{M}(E(A))$.*

PROOF. By the well-known duality between Boolean algebras and Boolean spaces [Halmos (1963), Theorem 6], $\mathfrak{M}(A) \approx \mathfrak{M}(\mathcal{B}_0(\mathfrak{M}(A)))$. Theorem 11.3 now yields the conclusion. ∎

**11.5 COROLLARY** *If $A$ is Boolean, then $A \cong \mathcal{B}_0(\mathfrak{M}(A))$; hence the completion of $A$ is $\mathcal{B}(\mathfrak{M}(A))$.*

The first assertion is the Stone representation theorem.

**11.6. ANNIHILATOR IDEALS.** Let $I$ be an ideal in a semi-prime ring $A$. We write

$$I\hat{\ } = \{a \in A : aI = 0\}.$$

An *annihilator ideal* is any ideal of the form $I\hat{\ }$. Clearly $I\hat{\ }\hat{\ } \supset I$, and $H \supset I$ implies $H\hat{\ } \subset I\hat{\ }$; it follows that an ideal $J$ is an annihilator ideal if and only if $J\hat{\ }\hat{\ } = J$. The family of all annihilator ideals in $A$ will be denoted by $\mathfrak{N}(A)$.

We define mappings $\Gamma$ and $\Delta$, as follows:

$$\Gamma I = \{M \in \mathfrak{M}(A) : M \not\supset I\} \qquad (I \text{ an ideal in } A);$$
$$\Delta \mathfrak{U} = \bigcap \mathfrak{U} \qquad (\mathfrak{U} \subset \mathfrak{M}(A)).$$

Since $\operatorname{cl} \mathfrak{U} = \{M : M \supset \bigcap \mathfrak{U}\}$,

$$(*) \qquad\qquad\qquad U\check{\ } = \Gamma \Delta \mathfrak{U}.$$

**11.7 LEMMA** [Lambek (1961), Theorem 6.6] *In a semi-prime ring $A$:*

(1) $\mathfrak{N}(A)$ *is a complete Boolean algebra, with set intersection as meet and $J \to J\hat{\ }$ as complementation.*

(2) $E(\mathfrak{Q}(A)) \cong \mathfrak{N}(\mathfrak{Q}(A))$, *under the mapping $e \to e\mathfrak{Q}(A)$.*

(3) $\mathfrak{N}(\mathfrak{Q}(A)) \cong \mathfrak{N}(A)$, *under the mapping $K \to K \cap A$.*

*If, in addition, $A$ is semi-simple, then*

(4) $I\hat{\ } = \Gamma \Delta I$.

11.8　Lemma　*If $A$ is semi-simple, then $\mathfrak{N}(A) \cong \mathcal{B}(\mathfrak{M}(A))$, under the mapping $\Gamma$.*

Proof. By (4) and (∗), $\Gamma I\hat{\ } = \Gamma\Delta\Gamma I = (\Gamma I)\check{\ } \in \mathcal{B}$; hence $\Gamma$ carries $\mathfrak{N}$ into $\mathcal{B}$. By (∗) and (4), $\Delta\mathfrak{U}\check{\ } = \Delta\Gamma\Delta\mathfrak{U} = (\Delta\mathfrak{U})\hat{\ } \in \mathfrak{N}$; so $\Delta$ carries $\mathcal{B}$ into $\mathfrak{N}$.

We now restrict $\Gamma$ to $\mathfrak{N}$ and $\Delta$ to $\mathcal{B}$. Clearly, $\Gamma(I \cap J) = \Gamma I \cap \Gamma J$; and as shown above, $\Gamma I\hat{\ } = (\Gamma I)\check{\ }$. Therefore $\Gamma$ is a homomorphism. In fact, $\Gamma$ is an isomorphism with inverse $\Delta\Gamma\Delta$. For, (∗) yields $\Gamma(\Delta\Gamma\Delta)\mathcal{B} = \mathcal{B}\check{\ }\check{\ } = \mathcal{B}$; and by (4), $(\Delta\Gamma\Delta)\Gamma J = J\hat{\ }\hat{\ } = J$. ∎

11.9　Theorem　*If $A$ is semi-simple, then $E(\mathfrak{Q}(A))$ is complete. In fact, $E(\mathfrak{Q}(A)) \cong \mathcal{B}(\mathfrak{M}(A))$, under the mapping*

$$e \to \Gamma(e\mathfrak{Q}(A) \cap A).$$

*Hence $E(\mathfrak{Q}(A))$ is the completion of $E(A)$.*

Proof. The isomorphism is obtained by combining (2) and (3) of 11.7 and 11.8. The concluding statement follows from 11.3. ∎

11.10　Lemma　*A space $X$ is extremally disconnected if and only if every regular open set in $X$ is closed.*

Proof. The usual definition of extremal disconnectedness is that the closure of each open set be open. (For its equivalence with the condition in 3.5, see [Gillman & Jerison (1960), 1H].) This implies that each regular open set $V$ is closed: $V = \operatorname{int}\operatorname{cl} V = \operatorname{cl} V$. Conversely, if $U$ is open, then $U\check{\ }$ is regular, so that if regular open sets are closed, then $\operatorname{cl} U$ is open. ∎

11.11　Corollary　*A necessary and sufficient condition that a semi-simple ring $A$ contain all idempotents of $\mathfrak{Q}(A)$ (i.e., that $E(A)$ be complete) is that $\mathfrak{M}(A)$ be extremally disconnected.*

Proof. The mapping $\Gamma$ of 11.6, when restricted to $E(A)$ (under the identification of $e$ with $(e)$), yields the mapping $\Gamma$ of 11.3. Hence by Theorem 11.9, $E(\mathfrak{Q}(A)) = E(A)$ if and only if $\mathcal{B}(\mathfrak{M}(A)) = \mathcal{B}_0(\mathfrak{M}(A))$, that is to say, every regular open set in $\mathfrak{M}(A)$ is closed. ∎

11.12　Corollary　*If $A$ is Boolean, then $\mathfrak{Q}(A) \cong \mathcal{B}(\mathfrak{M}(A))$, the completion of $A$.*

Proof. By Theorem 11.1, $\mathfrak{Q}(A) = E(\mathfrak{Q}(A))$; Theorem 11.9 now yields the isomorphsim and Corollary 11.5 the comment. ∎

It is shown in [Brainerd & Lambek (1959)] that the completion of a Boolean algebra $A$ is $\mathcal{Q}(A)$. By Corollary 11.5, this completion is $\mathcal{B}(\mathfrak{M}(A))$. Hence the present isomorphism provides a proof of either assuming the other.

**11.13  Corollary**  *For any space $X$,*

$$E(C^*(X)) \cong \mathcal{B}_0(X) \cong \mathcal{B}_0(\beta X)$$

*and*

$$E(\mathrm{Q}(X)) \cong \mathcal{B}(X) \cong \mathcal{B}(\beta X).$$

*Hence $C^*(X)$ contains all idempotents of $\mathrm{Q}(X)$ if and only if $X$ is extremally disconnected. (Cf. 4.12)*

**Proof.**  The assertions about $E(C^*)$ are obvious. (Cf. 11.3.) The rest then follow from 11.9 and (1) of 11.2. ∎

We remark that if $X$ is dense in $Y$, then $V \to V \cap X$ is one-one from $\mathcal{B}(Y)$ onto $\mathcal{B}(X)$, with inverse $U \to \mathrm{int}_Y \mathrm{cl}_Y U$. The mapping $V \to V \cap X$ carries $\mathcal{B}_0(Y)$ into $\mathcal{B}_0(X)$—and when $Y = \beta X$, onto.

**11.14  Corollary**  *If $A$ is semi-simple and contains $E(\mathcal{Q}(A))$, and if $\mathfrak{M}(A)$ is Hausdorff, then $\mathfrak{M}(A) \approx \mathfrak{M}(E(A))$.*

**Proof.**  By Corollary 11.11, the Hausdorff space $\mathfrak{M}(A)$ is extremally disconnected, hence totally disconnected. Apply 11.4. ∎

**11.15  Theorem**  *For any space $X$, $\mathfrak{M}(\mathcal{B}(X))$ is homeomorphic with the compact extremally disconnected space $K$ described in Theorem 6.9. Hence $\overline{\mathrm{Q}}^*(X) \cong C(\mathfrak{M}(\mathcal{B}(X)))$ and $\overline{\mathrm{Q}}(X) \cong \mathrm{Q}(\mathfrak{M}(\mathcal{B}(X)))$.*

**Proof.**  Applying 11.13, 10.17, 6.9 and the fact that $K$ is compact, we get, successively,

$$\mathfrak{M}(\mathcal{B}(X)) \approx \mathfrak{M}(E(\mathrm{Q}(X))) \approx \mathfrak{M}(\overline{\mathrm{Q}}^*(X)) \approx \mathfrak{M}(C(K)) \approx K.$$

The rest follows from 6.9. ∎

**11.16  Lemma**  *$\mathcal{B}(X)$ is the same for all (completely regular) spaces $X$ without isolated points and having a countable base. In fact, $\mathcal{B}(X) \cong \mathcal{Q}(B_\infty)$ where $B_\infty$ denotes the free Boolean algebra with countably infinitely many generators.*

**Proof.**  This is stated in [Birkhoff (1948), p. 177], except for two modifications. In [Birkhoff (1948)], $X$ is a $T_1$-space with a countable base of regular open sets; but a (completely) regular space with a countable base $\mathcal{U}$ does have a countable base of regular open sets—namely, $\{\mathrm{int}\,\mathrm{cl}\,U : U \in \mathcal{U}\}$. In [Birkhoff (1948)], our $\mathcal{Q}(B_\infty)$ is referred to as the completion of $B_\infty$; but by Corollary 11.12, the two are the same. ∎

11.17  THEOREM  $E(Q(X))^4$ *and* $\mathfrak{M}(Q(X))$ *are the same for all spaces* $X$ *without isolated points and having a countable base. In fact,* $E(Q(X)) \cong \mathcal{Q}(B_\infty)$ *and* $\mathfrak{M}(Q(X)) \approx \mathfrak{M}(\mathcal{Q}(B_\infty))$.

PROOF. By 11.13 and 11.16, $E(Q(X)) \cong \mathcal{Q}(B_\infty)$. And according to 10.14 (or by 11.14), $\mathfrak{M}(Q(X)) \approx \mathfrak{M}(E(Q(X)))$. ∎

---

[4]The original had $E(\mathcal{Q}(X))$.

# Index